\documentclass{amsart}
\usepackage{amsmath}
\usepackage{amssymb} 
\usepackage{amscd} 
\usepackage{graphicx} 
\usepackage{epsfig}
\usepackage[dvips]{color}
\usepackage{slashbox}

\newtheorem{theorem}{Theorem}[section]
\newtheorem{lemma}[theorem]{Lemma}

\newtheorem{proposition}[theorem]{Proposition}
\newtheorem{corollary}[theorem]{Corollary}
\newtheorem{claim}[theorem]{Claim}

\newtheorem{conjecture}[theorem]{Conjecture}
\newtheorem{remark}[theorem]{Remark}
\newtheorem{question}[theorem]{Question}

\theoremstyle{definition}
\newtheorem{definition}[theorem]{Definition}
\newtheorem{example}[theorem]{Example}

\numberwithin{equation}{section}
\numberwithin{figure}{section}
\numberwithin{table}{section}

\begin{document}
\baselineskip 15pt

\title{Left-orderable, non--$L$--space surgeries on knots}

\author[K. Motegi]{Kimihiko Motegi}
\address{Department of Mathematics, Nihon University, 
Tokyo 156--8550, Japan}
\email{motegi@math.chs.nihon-u.ac.jp}

\author[M. Teragaito]{Masakazu Teragaito}
\address{Department of Mathematics and Mathematics Education, Hiroshima University, 
1-1-1 Kagamiyama, Higashi-Hiroshima, 739--8524, Japan}
\email{teragai@hiroshima-u.ac.jp}

\thanks{
The first author has been partially supported by JSPS Grants--in--Aid for Scientific 
Research (C), 21540098, The Ministry of Education, Culture, Sports, Science and Technology, Japan and Joint Research Grant of Institute of Natural Sciences at Nihon University for 2013. \\
The second author has been partially supported by JSPS Grant-in-Aid for Scientific Research (C), 25400093, The Ministry of Education, Culture, Sports, Science and Technology, Japan.}

\dedicatory{In honor of Dale Rolfsen}

\begin{abstract}
Let $K$ be a knot in the $3$--sphere $S^3$. 
An $r$--surgery on $K$ is  \textit{left-orderable} 
if the resulting $3$--manifold $K(r)$ of the surgery has left-orderable fundamental group, 
and an $r$--surgery on $K$ is called an \textit{$L$--space surgery} 
if $K(r)$ is an $L$--space.  
A conjecture of Boyer, Gordon and Watson says that non-reducing surgeries on $K$ can be classified into left-orderable surgeries or $L$--space surgeries. 
We introduce a way to provide knots with left-orderable, non--$L$--space surgeries. 
As an application we present infinitely many hyperbolic knots on each of which 
every nontrivial surgery is a hyperbolic,  left-orderable, non--$L$--space surgery. 
\end{abstract}

\maketitle

{
\renewcommand{\thefootnote}{}
\footnotetext{
}

\section{Introduction}
\label{section:Introduction}

A nontrivial group $G$ is said to be \textit{left-orderable} if there exists a strict total ordering 
$<$ on its elements such that $g < h$ implies $fg < fh$ for all elements 
$f, g, h \in G$. 
The left-orderability of fundamental groups of $3$--manifolds has been studied by Boyer, Rolfsen and Wiest \cite{BRW}. 
In particular, 
they prove that the fundamental group of a $P^2$--irreducible $3$--manifold is left-orderable if and only if it has an epimorphism to a left-orderable group \cite[Theorem 1.1(1)]{BRW}. 
Since the infinite cyclic group $\mathbb{Z}$ is left-orderable, 
a $P^2$--irreducible $3$--manifold with first Betti number $b_1 \ge 1$ has left-orderable fundamental group. 
One obstruction for $G$ being left-orderable is an existence of torsion elements in $G$. 
Thus, for instance, 
lens spaces cannot have left-orderable fundamental groups. 
It is interesting to characterize rational homology $3$--spheres whose fundamental groups are left-orderable. 
Examples suggest that there exists a correspondence between rational homology $3$--spheres whose fundamental groups cannot be left-ordered and 
$L$--spaces which appear in the Heegaard Floer homology theory \cite{OS1, OS2}.  
For a rational homology $3$--sphere $M$, 
we have $\textrm{rk}\widehat{HF}(M) \ge  | H_1(M; \mathbb{Z}) |$. 
If the equality holds, 
then $M$ is called an \textit{$L$--space}. 
Following \cite[1.1]{BGW}, 
for homogeneity, we use $\mathbb{Z}_2$-coefficients for Heegaard Floer homology. 

The present paper is motivated by the following conjecture formulated by Boyer, Gordon and Watson \cite{BGW}. 

\begin{conjecture}[\cite{BGW}]
\label{LspaceConjecture}
An irreducible rational homology $3$--sphere is an $L$--space if and only if 
its fundamental group is not left-orderable. 
\end{conjecture}

In \cite{BGW} the conjecture is verified for geometric, non-hyperbolic $3$--manifolds 
and the $2$--fold branched covers of non-splitting alternating links. 
See also \cite{BB, CLW, Gr, I, P} for related results. 

A useful way to construct rational homology $3$--spheres is Dehn surgery on knots in the $3$--sphere $S^3$.  
Henceforth we will focus on Conjecture~\ref{LspaceConjecture} for rational homology $3$--spheres obtained by Dehn surgery on knots in $S^3$. 
For any knot $K$ in $S^3$
the exterior $E(K) = S^3 - \mathrm{int}N(K)$ has left-orderable fundamental group \cite[Corollary~3.5]{BRW}. 
However, 
the result $K(r)$ of $r$--Dehn surgery may not have such a fundamental group; 
see Examples~\ref{torus knots}, \cite{CW2} and \cite{Jun}. 

A Dehn surgery is said to be \textit{left-orderable} 
if the resulting $3$--manifold of the surgery has left-orderable fundamental group. 
Define the set of left-orderable surgeries on $K$ as    
$$\mathcal{S}_{LO}(K)= \{ r \in \mathbb{Q}\ |\ \pi_1(K(r))\ \textrm{is\ left-orderable} \}.$$ 
Similarly a Dehn surgery is called an \textit{$L$--space surgery} 
if the resulting $3$--manifold of the surgery is an $L$--space, 
and the set of $L$--space surgeries on $K$ is defined as   
$$\mathcal{S}_{L}(K)= \{ r \in \mathbb{Q}\ |\ K(r)\ \textrm{is\ an $L$--space} \}.$$ 

\begin{remark}
\label{mirror}
\begin{enumerate}
\item
Note that $0$--surgery does not yield a rational homology $3$--sphere, 
and hence $K(0)$ is not an $L$--space and $0 \not\in \mathcal{S}_{L}(K)$. 
On the other hand, 
if $K$ is a trivial knot, 
then $K(0) \cong S^2 \times S^1$ which has left orderable fundamental group. 
If $K$ is a nontrivial knot,  
then $K(0)$ is irreducible \cite[Corollary 8.3]{GabaiIII} and 
$H_1(K(0)) \cong \mathbb{Z}$, 
hence $0 \in  \mathcal{S}_{LO}(K)$ \cite[Theorem 1.1]{BRW}. 

\item
Let $K^*$ be the mirror image of a knot $K$, 
and put $-\mathcal{S} = \{ -r\ |\ r \in \mathcal{S}\}$ for $\mathcal{S} \subset \mathbb{Q}$. 
Since $K^*(-r)$ is orientation reversingly diffeomorphic to $K(r)$ and 
the conditions of a $3$--manifold $M$ having left-orderable fundamental group and 
being an $L$--space are independent of the orientation of $M$ \cite[p.1288]{OS3}, 
we have $\mathcal{S}_{LO}(K^*)= - \mathcal{S}_{LO}(K)$ and 
$\mathcal{S}_{L}(K^*)= - \mathcal{S}_{L}(K)$. 

\end {enumerate}
\end{remark}

If $K(r)$ is a reducible $3$--manifold for a nontrivial knot $K$, 
it has a lens space summand \cite[Theorem 3]{GL}, 
hence $r \not\in \mathcal{S}_{LO}(K)$, 
but $r$ may or may not be in $\mathcal{S}_L(K)$; 
see Remark~\ref{cable} and Example~\ref{torus knots}.  
 
If $K(r)$ is irreducible, 
Conjecture~\ref{LspaceConjecture} asserts that $r$ belongs to exactly one of 
$\mathcal{S}_{LO}(K)$ and $\mathcal{S}_L(K)$. 
Taking the cabling conjecture \cite{GS} into consideration, Conjecture~\ref{LspaceConjecture} suggests: 

\begin{conjecture}
\label{LspaceConjecture2}
Let $K$ be a knot in $S^3$ which is not a cable of a nontrivial knot. 
Then 
$\mathcal{S}_{LO}(K) \cup \mathcal{S}_{L}(K) = \mathbb{Q}$ and 
$\mathcal{S}_{LO}(K) \cap \mathcal{S}_{L}(K) = \emptyset$. 
\end{conjecture}

\begin{remark}
\label{cable}
The cabling conjecture \cite{GS} asserts that if $K(r)$ is reducible for a nontrivial knot $K$, 
then $K$ is cabled and $r$ is a cabling slope. 
There exists a cable knot $K$ 
for which $\mathcal{S}_{LO}(K) \cup \mathcal{S}_{L}(K) \ne \mathbb{Q}$. 
For instance, 
let $K$ be a $(p, q)$ cable of a non-fibered knot $k$ $(q > 0)$. 
Then $K(pq) = k(\frac{p}{q}) \sharp L(q, p)$ \cite[Corollary 7.3]{Go}. 
Since $\pi_1(K(pq))$ has a torsion, 
$pq \not\in \mathcal{S}_{LO}(K)$. 
To see that $pq \not\in \mathcal{S}_{L}(K)$, 
we note that  $\widehat{HF}(K(pq)) \cong \widehat{HF}(k(\frac{p}{q})) \otimes \widehat{HF}(L(q, p))$; 
see \cite[8.1(5)]{Szabo} $($\cite{OS2}$)$.  
Since $k$ is a non-fibered knot, 
$k(\frac{p}{q})$ is not an $L$--space \cite{Ni, Ni2}. 
Hence the rank of $\widehat{HF}(K(pq))$ is strictly bigger than $|p|q$, 
and $K(pq)$ is not an $L$--space. 
It follows that $pq \not\in \mathcal{S}_{LO}(K) \cup \mathcal{S}_{L}(K)$. 
\end{remark}

\medskip

For the trivial knot and nontrivial torus knots, 
Examples~\ref{trivial knot} and \ref{torus knots} describe 
$\mathcal{S}_{LO}(K)$ and $\mathcal{S}_{L}(K)$ explicitly. 
Note that these knots satisfy Conjecture~\ref{LspaceConjecture2}. 

\begin{example}[\textbf{trivial knot}]
\label{trivial knot}
Let $K$ be the trivial knot in $S^3$. 
Then $\mathcal{S}_{LO}(K) = \{ 0 \}$ and 
$\mathcal{S}_{L}(K) =\mathbb{Q} - \{ 0 \}$.
\end{example}

\begin{example}[\textbf{torus knots}]
\label{torus knots}
For a nontrivial torus knot $T_{p, q}$ $(p > q \ge 2)$, 
the argument in the proof of \cite[Theorem 1.4]{CW1} shows that  
$\mathcal{S}_{LO}(T_{p, q})= (-\infty,\ pq-p-q) \cap \mathbb{Q}$ and 
$\mathcal{S}_{L}(T_{p, q})=[pq-p-q,\ \infty) \cap \mathbb{Q}$. 
\end{example}

\begin{example}[\textbf{figure-eight knot}]
\label{figure-eight knot}
Let $K$ be the figure-eight knot. 
Following \cite{OS3, OS4}, 
$\mathcal{S}_{L}(K)=\emptyset$. 
Thus it is expected that 
$\mathcal{S}_{LO}(K)=\mathbb{Q}$. 
Boyer, Gordon and Watson \cite{BGW} show that 
$\mathcal{S}_{LO}(K) \supset (-4, 4) \cap \mathbb{Q}$, 
and Clay, Lidman and Watson \cite{CLW} improve that $\mathcal{S}_{LO}(K) \supset [-4, 4] \cap \mathbb{Q}$. 
Furthermore, \cite{F} implies that $\mathcal{S}_{LO}(K) \supset \mathbb{Z}$.  
\end{example}

For related results, 
see  \cite{CT, CW2, HT2, LR, Tera, Tran2}. 

It is known that there exist some constraints for knots which admit $L$--space surgeries. 
For instance, such knots have specific Alexander polynomials \cite{OS3}, and 
must be fibered \cite{Ni, Ni2}. 
Thus generically we have $\mathcal{S}_L(K) = \emptyset$. 
Hence Conjecture~\ref{LspaceConjecture2} suggests that 
$\mathcal{S}_{LO}(K) = \mathbb{Q}$ for most knots. 
Despite being expected, 
there is no literature giving explicitly knots with $\mathcal{S}_{LO}(K) = \mathbb{Q}$ and 
$\mathcal{S}_L(K) = \emptyset$. 
In the present note we give infinitely many 
satellite knots and hyperbolic knots with this property.  

\begin{theorem}
\label{satellite}
Given a nontrivial knot $K'$, 
there are infinitely many prime satellite knots $K$ each of which has $K'$ as a companion knot and enjoys the following properties:  
\begin{enumerate}
\item
$K(r)$ is a toroidal $3$--manifold which is not a graph manifold for  all but finitely many $r \in \mathbb{Q}$. 
\item
$\mathcal{S}_{LO}(K) = \mathbb{Q}$. 
\item
$\mathcal{S}_L(K) = \emptyset$.
\end{enumerate}
\end{theorem}

This is an application of Proposition~\ref{ClayWatson} due to Clay and Watson 
\cite[Proposition~4.1]{CW1}. 
In Theorem~\ref{satellite}, 
$K$ is satellite knot and the resulting $3$--manifold $K(r)$ has a nontrivial Jaco-Shalen-Johannson (JSJ) decomposition \cite{JS, Jo}. 
Since Proposition~\ref{ClayWatson} does not work for creation of hyperbolic knots, 
we will introduce an 
effective way to provide infinitely many hyperbolic knots having left-orderable, non--$L$--space surgeries from a given knot with 
left-orderable surgeries; 
see Section~\ref{constructions}. 
Then we will apply the construction to prove the following: 

\begin{theorem}
\label{all rationals}
There exist infinitely many hyperbolic knots $K$ each of which enjoys the following properties. 
\begin{enumerate}
\item 
$K(r)$ is a hyperbolic $3$--manifold for all $r \in \mathbb{Q}$. 
\item 
$\mathcal{S}_{LO}(K) =  \mathbb{Q}$. 
\item
$\mathcal{S}_{L}(K) = \emptyset$. 
\end{enumerate}
\end{theorem}

\bigskip

\section{Left-orderable surgeries on periodic knots}
\label{LO_periodic}

A knot $K$ in $S^3$ is called a \textit{periodic knot} with period $p$ 
if there is an orientation preserving diffeomorphism 
$f : S^3 \to S^3$ such that $f(K) = K$, $f^p = id\ (p > 1)$, $\mathrm{Fix}(f) \ne \emptyset$, 
and $\mathrm{Fix}(f) \cap K = \emptyset$, 
where  $\mathrm{Fix}(f)$ is the set of fixed points of $f$. 
By the positive answer to the Smith conjecture \cite{MB}, 
$f$ is a rotation of $S^3$ about the unknotted circle $C = \mathrm{Fix}(f)$. 
So by taking the quotient $S^3/ \langle f \rangle$, 
we obtain the \textit{factor knot} $\overline{K} = K/ \langle f \rangle$ 
and the unknotted circle $\overline{C} = C/\langle f \rangle$  in $S^3 = S^3 / \langle f \rangle$. 
We often call $C$ the \textit{axis} and $\overline{C}$ the \textit{branch circle}. 
Since $K$ is connected, 
the linking number $lk(\overline{K}, \overline{C})$ and the period $p$ are relatively prime. 
Note that if the periodic knot $K$ is unknotted, 
then the equivariant loop theorem \cite{MY} implies that 
$K \cup C$ is the Hopf link and $\overline{K} \cup \overline{C}$ is also the Hopf link. 
To exclude such a trivial case, 
in the following we consider nontrivial periodic knots. 

\begin{figure}[htbp]
\begin{center}
\includegraphics[width=0.65\linewidth]{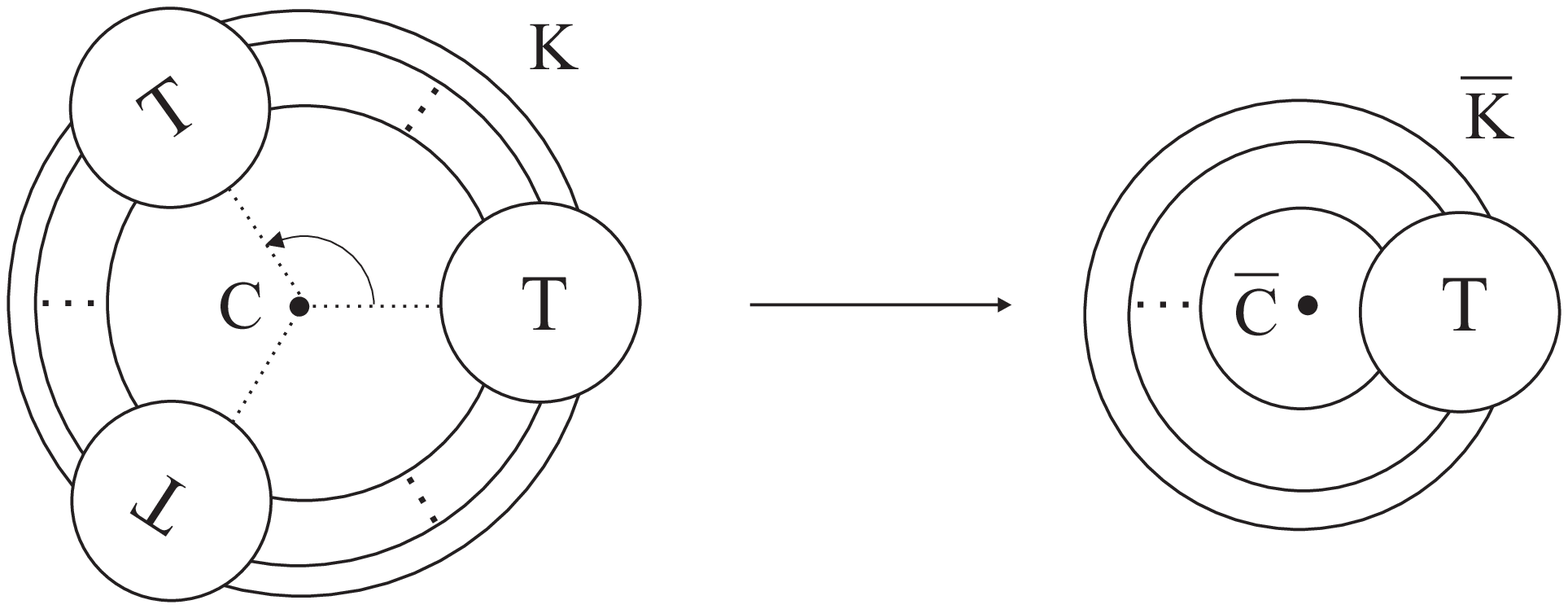}
\caption{A periodic knot $K$ with an axis $C$ and its factor knot $\overline{K}$; $T$ is a tangle.}
\label{factorknot1}
\end{center}
\end{figure}

The next theorem asserts that 
the fundamental groups of $3$--manifolds obtained by Dehn surgeries on the periodic knot $K$ inherit left-orderability from those of $3$--manifolds obtained by Dehn surgeries on the factor knot $\overline{K}$. 
For a subset $\mathcal{S} \subset \mathbb{Q}$ and 
a positive integer $p$, 
we denote by $p\mathcal{S}$ the subset 
$\{ pr\ |\ r \in \mathcal{S} \} \subset \mathbb{Q}$. 
Note that if $\mathcal{S} = \mathbb{Q}$, 
then $p \mathcal{S} = \mathbb{Q}$.  

\begin{theorem}
\label{factorknot_LO}
Let $K$ be a nontrivial knot in $S^3$ with cyclic period $p$, 
and let $\overline{K}$ be its factor knot.  
Then $\mathcal{S}_{LO}(K) \supset p\mathcal{S}_{LO}(\overline{K})$. 
\end{theorem}

\noindent
\textit{Proof of Theorem~\ref{factorknot_LO}.}
Let $f : S^3 \to S^3$ be an orientation preserving diffeomorphism giving the cyclic period $p$ of $K$ 
with axis $C = \mathrm{Fix}(f)$ and factor knot $\overline{K} = K / \langle f \rangle$. 
Take an $\langle f \rangle$--invariant tubular neighborhood $N(K)$ of $K$. 
Let $N(\overline{K})$ be the quotient $N(K)/ \langle f \rangle$. 
In the following $E(K) = S^3 - \mathrm{int}N(K)$ and $E(\overline{K}) = E(K) / \langle f \rangle = S^3 - \mathrm{int}N(\overline{K})$. 
Denote by $(\mu, \lambda)$ (resp. $(\overline{\mu}, \overline{\lambda})$) 
a preferred meridian-longitude pair of $\pi_1(\partial N(K))$ (resp. $\pi_1(\partial N(\overline{K}))$).  
We can choose a simple closed curve representing the preferred longitude 
$\lambda$ which is invariant under $\langle f \rangle$; see \cite{E}. 

Let $\pi : E(K) \to E(\overline{K})$ be the cyclic branched covering branched along $\overline{C} = C / \langle f \rangle$. 

\begin{lemma}
\label{extension}
The branched cover $\pi : E(K) \to E(\overline{K})$ can be extended to 
a branched cover $\pi' : K(\frac{m}{n}) \to \overline{K}(\frac{m}{pn})$. 
\end{lemma}

\noindent
\textit{Proof of Lemma~\ref{extension}.}
Let $\pi_* : \pi_1(E(K)) \to \pi_1(E(\overline{K}))$ be the homomorphism induced by $\pi$. 
Then $\pi_*|_{\pi_1(\partial E(K))} : \pi_1(\partial E(K)) \to \pi_1(\partial E(\overline{K}))$ 
sends $\mu$ to $\overline{\mu}$, and $\lambda$ to $p \overline{\lambda}$. 
Hence 
$\pi_*|_{\pi_1(\partial E(K))}(m \mu + n \lambda) = m \overline{\mu} + pn \overline{\lambda} 
= (m,\ p)\{ \frac{m}{(m,\ p)} \overline{\mu} + \frac{pn}{(m,\ p)}\overline{\lambda} \}$. 
Then we can extend $\pi : E(K) \to E(\overline{K})$ to 
${\pi}' : K(\frac{m}{n}) = E(K) \cup V \to \overline{K}(\frac{m}{pn}) = E(\overline{K}) \cup \overline{V}$, 
where $V,\ \overline{V}$ are filled solid tori. 
If $(m, p) \ge 2$, then ${\pi}'$ branches along the core of the filled solid torus $\overline{V}$ as well as $\overline{C}$. 
\hspace*{\fill} $\square$(Lemma~\ref{extension})

Thus we have a commutative diagram: 

\begin{eqnarray*}
\begin{CD}
	E(K) @>  \pi >>E(\overline{K}) \\
@V \mathrm{Dehn\ filling}VV
		@VV  \mathrm{Dehn\ filling} V\\
	K(\frac{m}{n}) @>>{\pi}'> \overline{K}(\frac{m}{pn})
\end{CD}
\end{eqnarray*}

\medskip

Assume that $\frac{m}{pn} \in \mathcal{S}_{LO}(\overline{K})$, 
i.e. 
$\overline{K}(\frac{m}{pn})$ has left-orderable fundamental group. 
Let us prove that $K(\frac{m}{n})$ has also left-orderable fundamental group, 
i.e. $p \times \frac{m}{pn} = \frac{m}{n} \in \mathcal{S}_{LO}(K)$.  

\begin{lemma}
\label{irreducible}
$K(\frac{m}{n})$ is irreducible.
\end{lemma}

\noindent
\textit{Proof of Lemma~\ref{irreducible}.}
Suppose for a contradiction that $K(\frac{m}{n})$ is reducible.
Since $K$ is a nontrivial periodic knot,  
$K$ is cabled and $\frac{m}{n}$ is the cabling slope \cite{LZ, HM, HS}. 
First we assume that $K$ is a torus knot. 
Then $E(K)$ has a unique Seifert fibration (up to isotopy). 
Following \cite[Theorem 2.2]{MS}, 
we choose a Seifert fibration of $E(K)$ which is preserved by $f$. 
If $C$ is not a fiber, 
we take a regular fiber $t$ intersecting $C$. 
Since $f$ fixes a point in $t \cap C$, 
$f(t) = t$ and $f$ reverses the orientation of $t$. 
This then implies that $f$ reverses the orientation of $K$, 
and hence $C$ intersects $K$ in exactly two points, a contradiction. 
Thus $C$ is a fiber in the $\langle f \rangle$--invariant Seifert fibration of $E(K)$. 
Since a regular fiber is knotted in $S^3$, 
$C$ is one of two exceptional fibers in $E(K)$. 
Then the quotient $E(\overline{K}) = E(K) / \langle f \rangle$ 
has also a Seifert fibration induced from that of $E(K)$ and thus $\overline{K}$ is a torus knot; 
the surgery slope $\frac{m}{pn}$ on $\partial E(\overline{K})$ is the fiber slope. 
Since $\frac{m}{pn}$ is the fiber (i.e. cabling) slope, 
$\pi_1(\overline{K}(\frac{m}{pn}))$ has a nontrivial torsion, 
contradicting the left-orderability of $\pi_1(\overline{K}(\frac{m}{pn}))$. 

Next assume that $K$ is an $(x, y)$--cable in a knotted solid torus $W$, where $y \ge 2$.   
By the $\langle f \rangle$--invariant version (\cite[Theorem 8.6]{MS}) of the torus decomposition theorem \cite{JS, Jo},  
we may assume that $f$ leaves a companion solid torus $W$ invariant. 
First we note that $W \cap C = \emptyset$. 
For otherwise, 
$f |_{\partial W}$ has fixed points and hence it is an involution, 
and $f$ reverses the orientation of an $\langle f \rangle$--invariant core of $W$.  
Hence it also reverses the orientation of $K$ (which has winding number $y \ge 2$ in $W$). 
This then implies that $C$ intersects $K$ in exactly two points, a contradiction. 
Thus $W \subset S^3 - C$.   
We denote the quotient $W / \langle f \rangle$ by $\overline{W}$. 
We may assume that the cable space $W - \mathrm{int}N(K)$ has a Seifert fibration preserved by $f$ \cite[Theorem 2.2]{MS}. 
Then $\overline{W} - \mathrm{int}N(\overline{K}) = (W - \mathrm{int}N(K)) / \langle f \rangle$ 
has an induced Seifert fibration in which a regular fiber on $\partial N(\overline{K})$ 
represents the surgery slope $\frac{m}{pn}$. 
This implies that the result of $\frac{m}{pn}$--surgery of $\overline{W}$ along $\overline{K}$, 
and hence $\overline{K}(\frac{m}{pn})$, 
has a nontrivial lens space summand whose fundamental group has order $y \ge 2$. 
Since $\pi_1(\overline{K}(\frac{m}{pn}))$ has a nontrivial torsion,  
it cannot be left-orderable, contradicting the assumption.  
\hspace*{\fill} $\square$(Lemma~\ref{irreducible})

\medskip

The above diagram induces the commutative diagram of fundamental groups below. 

\begin{eqnarray*}
\begin{CD}
	\pi_1(E(K)) @>  \pi_* >>\pi_1(E(\overline{K})) \\
@V VV
		@VV  V\\
	\pi_1(K(\frac{m}{n})) @>>{\pi}'_*> \pi_1(\overline{K}(\frac{m}{pn}))
\end{CD}
\end{eqnarray*}

\smallskip

\begin{lemma}
\label{surjective}
$\pi'_* : \pi_1(K(\frac{m}{n})) \to \pi_1(\overline{K}(\frac{m}{pn}))$ is surjective. 
\end{lemma}

\noindent
\textit{Proof of Lemma~\ref{surjective}.}
Choose a point $x \in C = \mathrm{Fix}(f)$ (resp. $\pi(x) \in \overline{C}$) 
as a base point of $\pi_1(E(K))$ (resp. $\pi_1(E(\overline{K})$). 
We take obvious meridians $\overline{\mu}_i$ of $\overline{K}$ which are generators of $\pi_1(E(\overline{K}), \pi(x))$ 
(with respect to the Wirtinger presentation of $\pi_1(E(\overline{K}), \pi(x))$). 
Then their lifts $\mu_i \in \pi_1(E(K))$ satisfy $\pi_*(\mu_i) = \overline{\mu}_i$, 
and hence $\pi_* : \pi_1(E(K)) \to \pi_1(E(\overline{K}))$ is an epimorphism. 
Since vertical homomorphisms are also epimorphisms, 
${\pi}'_* : \pi_1(K(\frac{m}{n})) \to \pi_1(\overline{K}(\frac{m}{pn}))$ is also an epimorphism. 
\hspace*{\fill} $\square$(Lemma~\ref{surjective})

\medskip

By Lemma~\ref{irreducible} $K(\frac{m}{n})$ is irreducible, 
and by Lemma~\ref{surjective} we have an epimorphism 
from $\pi_1(K(\frac{m}{n}))$ to the left-orderable group $\pi_1(\overline{K}(\frac{m}{pn}))$.  
Then it follows from \cite[Theorem 1.1(1)]{BRW} that $\pi_1(K(\frac{m}{n}))$ is also left-orderable.  
Thus if $r = \frac{m}{pn} \in \mathcal{S}_{LO}(\overline{K})$, 
then $pr = \frac{m}{n} \in \mathcal{S}_{LO}(K)$. 
\hspace*{\fill} $\square$(Theorem~\ref{factorknot_LO})

\bigskip

\section{$L$--space surgeries on periodic knots}
\label{Lspaces}

In \cite{Ni, Ni2} Ni proves that if a knot $K$ in $S^3$ has an $L$--space surgery, 
then $K$ is a fibered knot, 
i.e. $E(K)$ has a fibering over the circle.  
For a periodic knot $K$, 
the next theorem gives a necessary condition on the factor knot for $K$ having 
an $L$--space surgery.   

\begin{theorem}
\label{nonLspaces}
Let $K$ be a periodic knot in $S^3$ with axis $C$, 
and let $\overline{K}$ be its factor knot with branch circle $\overline{C}$. 
Suppose that $K$ has an $L$--space surgery. 
Then $E(\overline{K})$ has a fibering over the circle with a fiber surface $\overline{S}$ such that 
$| \overline{S} \cap \overline{C} |$ equals the algebraic intersection number between $\overline{S}$ and $\overline{C}$, 
i.e. the linking number $lk(\overline{K}, \overline{C})$. 
\end{theorem}

In particular, we have: 

\begin{corollary}
\label{factor_nonfibered}
Let $K$ be a periodic knot with factor knot $\overline{K}$. 
If $\overline{K}$ is not fibered, 
then $\mathcal{S}_L(K) = \emptyset$. 
\end{corollary}

\noindent
\textit{Proof of Theorem~\ref{nonLspaces}.}
Let $f : S^3 \to S^3$ be an orientation preserving diffeomorphism of finite order satisfying $f(K) =K$. 
Note that $C = \mathrm{Fix}(f)$, $\overline{K} = K / \langle f \rangle$ and $\overline{C} = C / \langle f \rangle$. 
Let $N(K)$ be an $\langle f \rangle$--invariant tubular neighborhood of $K$. 

Assume that $K$ has an $L$--space surgery. 
Then Ni \cite[Corollary 1.3]{Ni} (\cite{Ni2}) proves that $E(K) = S^3 - \mathrm{int}N(K)$ 
has a fibering over the circle. 
Following Proposition 6.1 in \cite{EL}, 
we can isotope the fibering to a fibering preserved by 
the action of $\langle f \rangle$ so that the axis $C$ is transverse to the fibers. 
Thus $E(\overline{K})$ inherits a fibering over the circle 
such that all the fibers are transverse to the branch circle $\overline{C} = C/ \langle f \rangle$. 
Let $\overline{S}$ be a fiber surface of $E(\overline{K})$. 
Since $\overline{C}$ intersects each fiber surface of the fibering of $E(\overline{K})$ transversely, 
$| \overline{S} \cap \overline{C} |$ coincides with the algebraic 
intersection number between $\overline{S}$ and $\overline{C}$, 
i.e. the linking number $lk(\partial \overline{S}, \overline{C})$, 
which equals the linking number $lk(\overline{K}, \overline{C})$. 
\hspace*{\fill} $\square$(Theorem~\ref{nonLspaces})

\medskip

As Ni \cite{Ni, Ni2} proves, 
the fiberedness of $K$ is necessary for $K$ to have an $L$--space surgery. 
On the other hand, 
the periodicity of $K$ itself also puts strong restrictions on $3$--manifolds obtained by 
Dehn surgeries on $K$. 
For instance, 
if a periodic knot $K$ with period $p > 2$ has a finite surgery, 
which is also an $L$--space surgery, 
then $K$ is a torus knot or a cable of a torus knot \cite[Proposition 5.6]{MM3}. 
So we would like to ask: 

\begin{question}
Let $K$ be a knot in $S^3$ with cyclic period $p > 2$ other than a torus knot or a cable of a torus knot. 
Then does $K$ admit an $L$--space surgery? 
\end{question}

\bigskip

\section{Periodic constructions}
\label{constructions}

Given a periodic knot, 
taking the quotient by the periodic automorphism, 
we obtain its factor knot; see Section~\ref{LO_periodic}. 
Reversing this procedure, 
we have: 

\begin{definition}[\textbf{periodic construction}]
\label{periodic construction}
Let $( \overline{K}, \overline{C} )$ be a pair of a knot $\overline{K}$ and an unknotted circle $\overline{C}$ 
which is disjoint from $\overline{K}$. 
Then for an integer $p \ge 2$ with $(p, lk(\overline{K}, \overline{C})) = 1$,  
take the $p$--fold cyclic branched cover of $S^3$ branched along $\overline{C}$ to 
obtain a knot $K_{\overline{C}}^p$ which is the preimage of $\overline{K}$. 
We call $K_{\overline{C}}^p$ the knot obtained from the pair $( \overline{K}, \overline{C} )$ by \textit{$p$--periodic construction}. 
\end{definition}

\medskip

Note that $K_{\overline{C}}^p$ is a knot with cyclic period $p$ whose factor knot is $\overline{K}$. 
Hence Theorems~\ref{factorknot_LO} and \ref{nonLspaces} immediately imply the following result. 

\begin{theorem}
\label{construction}
Let $( \overline{K}, \overline{C} )$ be a pair as in Definition~\ref{periodic construction}. 
If $\overline{K}$ is a fibered knot, 
$\overline{C}$ is chosen so that any fiber surface $($i.e. minimal genus Seifert surface$)$ $\overline{S}$ satisfies 
the inequality $| \overline{S} \cap \overline{C} | > lk(\overline{K}, \overline{C})$. 
Then a knot $K_{\overline{C}}^p$ obtained from the pair $( \overline{K}, \overline{C} )$ by $p$--periodic construction 
enjoys the following properties: 

\begin{enumerate}
\item
$\mathcal{S}_{LO}(K_{\overline{C}}^p) \supset p \mathcal{S}_{LO}(\overline{K})$. 
\item
$\mathcal{S}_{L}(K_{\overline{C}}^p) = \emptyset$. 
\end{enumerate}
\end{theorem}

\medskip

If $\overline{K}$ is a trivial knot, 
then $\mathcal{S}_{LO}(\overline{K}) = \{ 0 \}$ 
and hence $p \mathcal{S}_{LO}(\overline{K}) = \{ 0 \}$. 
So we will apply Theorem~\ref{construction} to nontrivial knots. 

\begin{remark}
\label{genus}
We denote the genus of a knot $k$ in $S^3$ by $g(k)$. 
For $\overline{K}$ and $K_{\overline{C}}^p$, 
we have $g(K_{\overline{C}}^p) \ge pg(\overline{K})$ \cite[Theorem~3.2]{Naik}. 
In particular, 
for a nontrivial knot $\bar{K}$, 
$g(K_{\overline{C}}^p) \to \infty$ as $p \to \infty$. 
\end{remark}

Theorem~\ref{construction} is accompanied by the following theorems. 

\begin{theorem}
\label{hyperbolic axis1}
Given a nontrivial knot $\overline{K}$ in $S^3$, 
we can take an unknotted circle $\overline{C}$ so that $\overline{K} \cup \overline{C}$ is a hyperbolic link
with arbitrary linking number.  
\end{theorem}

\noindent
\textit{Proof of Theorem~\ref{hyperbolic axis1}.}
The following argument is based on the proofs of Theorems 2.1 and 2.2 in \cite{Ai}.    
Arrange $\overline{K}$ as a closed $n$--braid for some integer $n$. 
If necessary, stabilizing the braid, 
we may assume that the braid contains both a positive crossing and a negative crossing (Figure~\ref{altrenatingcrossings}). 
Then introduce $(n-1)$--strands $C_i$ $(i = 1, \dots, n-1)$ 
between the $n$--strands of the original braid so that the crossings introduced, 
together with the original crossings, 
are alternately positive and negative. 
See Figure~\ref{altrenatingcrossings}. 

\begin{figure}[htbp]
\begin{center}
\includegraphics[width=0.7\linewidth]{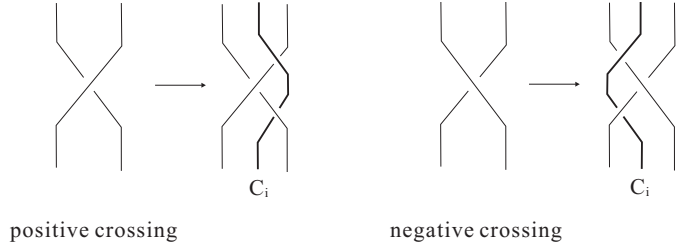}
\caption{Insertion of $(n-1)$--strands; $n = 2$}
\label{altrenatingcrossings}
\end{center}
\end{figure}

Then we arrange $C_i$ as in Figure~\ref{connect} 
so that the closed braid is a $2$--component link consisting of $\overline{K}$ and 
an unknotted circle $\overline{C} = C_1 \cup \cdots \cup C_{n-1}$ and $\overline{K} \cup \overline{C}$ 
is a non-split prime alternating link \cite[Theorem 1]{Me}. 

\begin{figure}[htbp]
\begin{center}
\includegraphics[width=0.85\linewidth]{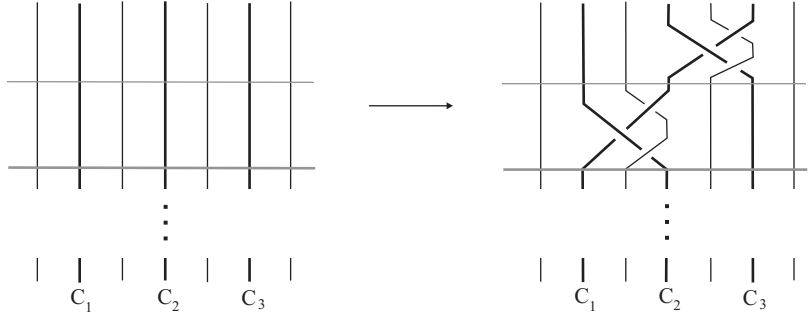}
\caption{Arrangement of $C_1, \dots, C_{n-1}$; $n=4$}
\label{connect}
\end{center}
\end{figure}

Since our braid contains both a positive crossing and a negative crossing, 
we can add some negative twists or  positive twists as in Figure~\ref{alternatingtwistings} to make $\overline{C}$ so that 
$lk(\overline{K}, \overline{C}) = l$ for a given integer $l$. 

\begin{figure}[htbp]
\begin{center}
\includegraphics[width=0.75\linewidth]{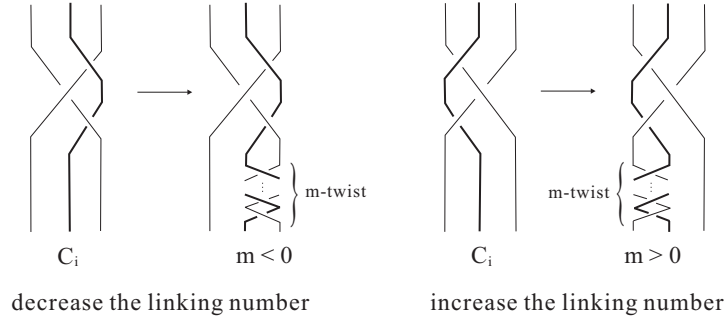}
\caption{Adding negative or positive twistings}
\label{alternatingtwistings}
\end{center}
\end{figure}

Note that the resulting link $\overline{K} \cup \overline{C}$ is also a non-split prime alternating link. 
It follows from \cite[Corollary 2]{Me} that $\overline{K} \cup \overline{C}$ is either a torus link or a hyperbolic link. 
Since $\overline{K}$ is nontrivial, but $\overline{C}$ is trivial, 
the former cannot occur, 
and thus $\overline{K} \cup \overline{C}$ is a hyperbolic link. 
\hspace*{\fill} $\square$(Theorem~\ref{hyperbolic axis1})

\begin{theorem}
\label{hyperbolic axis2}
\begin{enumerate}
\item
If $\overline{K} \cup \overline{C}$ is a hyperbolic link and $p > 2$, 
then $K_{\overline{C}}^p$ is a hyperbolic knot, 
and $K_{\overline{C}}^p(r)$ is a hyperbolic $3$--manifold for all $r \in \mathbb{Q}$. 
\item
Assume that $p > 2$ and $\overline{C_i}$ $(i = 1, 2)$ is an unknotted circle such that 
$lk(\overline{K}, \overline{C_i})$ and $p$ are relatively prime, 
and $\overline{K} \cup \overline{C_i}$ is a hyperbolic link.  
If $K_{\overline{C_1}}^p$ and $K_{\overline{C_2}}^p$ are isotopic in $S^3$, 
then $\overline{K} \cup \overline{C_1}$ and $\overline{K} \cup \overline{C_2}$ are isotopic.  
\end{enumerate}
\end{theorem}

\noindent
\textit{Proof of Theorem~\ref{hyperbolic axis2}.}
$(1)$ Assume for a contradiction that $K_{\overline{C}}^p$ is not hyperbolic. 
Then it is either a torus knot or a satellite knot. 
Let $f : S^3 \to S^3$ be the deck transformation of the $p$--fold cyclic branched cover given in 
Theorem~\ref{construction}, 
which is an orientation preserving diffeomorphism giving the cyclic period $p$ of $K_{\overline{C}}^p$. 
In the following, 
we take an $\langle f \rangle$--invariant tubular neighborhood $N(K_{\overline{C}}^p)$ and 
denote $S^3 - \mathrm{int}N(K_{\overline{C}}^p)$ by $E(K_{\overline{C}}^p)$.  
The preimage of the branch circle $\overline{C}$ is an unknotted circle $C = \mathrm{Fix}(f)$, 
which is contained in the interior of $E(K_{\overline{C}}^p)$.  
Note also that $K_{\overline{C}}^p$ is a nontrivial knot. 
For otherwise, the equivariant loop theorem \cite{MY} implies that $K_{\overline{C}}^p \cup C$ is the Hopf link and $\overline{K} \cup \overline{C}$ is also the Hopf link, 
contradicting the hyperbolicity of $\overline{K} \cup \overline{C}$. 

\medskip

\begin{claim}
\label{nontorus}
$K_{\overline{C}}^p$ is not a torus knot. 
\end{claim}

\noindent
\textit{Proof of Claim~\ref{nontorus}.}
Assume for a contradiction that $K_{\overline{C}}^p$ is a torus knot. 
Then $E(K_{\overline{C}}^p)$ has a unique Seifert fibration up to isotopy. 
We choose a Seifert fibration of $E(K_{\overline{C}}^p)$ which is preserved by $f$ \cite[Theorem 2.2]{MS}.  
Then the argument in the proof of Lemma~\ref{irreducible} shows that 
$C$ is one of two exceptional fibers in $E(K_{\overline{C}}^p)$. 
Then the quotient 
$E(\overline{K}) - \mathrm{int}N(\overline{C}) 
= (E(K_{\overline{C}}^p) - \mathrm{int}N(C)) / \langle f \rangle$ 
has also a Seifert fibration. 
Thus $S^3 - \mathrm{int}N(\overline{K} \cup \overline{C}) = E(\overline{K}) - \mathrm{int}N(\overline{C})$ is a Seifert fiber space, 
contradicting its hyperbolicity.  
\hspace*{\fill} $\square$(Claim~\ref{nontorus})

\medskip

\begin{claim}
\label{nonsatellite}
$K_{\overline{C}}^p$ is not a satellite knot. 
\end{claim}

\noindent
\textit{Proof of Claim~\ref{nonsatellite}.}
Suppose for a contradiction that $K_{\overline{C}}^p$ is a satellite knot. 
Then we have an $\langle f \rangle$--invariant torus decomposition of $E(K_{\overline{C}}^p)$ \cite[Theorem 8.6]{MS}. 
Let $\Sigma$ be the invariant family of essential tori in $E(K_{\overline{C}}^p)$. 

\textit{Case $(i)$.}
There is an essential torus $T \in \Sigma$ such that $f(T) = T$. 
Then $T$ bounds an $\langle f \rangle$--invariant companion solid torus $W$ containing $K_{\overline{C}}^p$. 
Note that $K_{\overline{C}}^p$ is not a core of $W$. 
We see that $W \cap C = \emptyset$, 
for otherwise $f|_{\partial W}$ has a fixed point and it is an involution, 
i.e. $(f|_{\partial W})^2$ is the identity map. 
By the classical Smith theory \cite{Smith} $f$ itself is an involution, 
contradicting the assumption. 
Thus $W$ lies in $S^3 - C$. 
We may assume that $W \subset S^3 - \mathrm{int}N(C)$ for a small tubular neighborhood 
$N(C)$ of $C$. 
Since the core of $W$ is not a core of $S^3 - \mathrm{int}N(C)$, 
$S^3 - \mathrm{int}N(K_{\overline{C}}^p \cup C)$ contains the $\langle f \rangle$--invariant essential torus $T = \partial W$. 
This then implies that  $S^3 - \mathrm{int}N(\overline{K} \cup \overline{C})$ contains an essential torus $\partial W/ \langle f \rangle$. 
This contradicts the hyperbolicity of $S^3 - \mathrm{int}N(\overline{K} \cup \overline{C})$. 

\textit{Case $(ii)$.}
For each $T \in \Sigma$, 
$f(T) \ne T$ (hence, $f(T) \cap T = \emptyset$). 
Let us pick an essential torus $T \in \mathcal{T}$. 
Note that $T$ is essential in $S^3 - \mathrm{int}N(K_{\overline{C}}^p \cup C)$. 
Then the image $\overline{T} \subset E(\overline{K} \cup \overline{C})$ of $T$ by the covering projection 
is also essential. 
This contradicts the hyperbolicity of $S^3 - \mathrm{int}N(\overline{K} \cup \overline{C})$. 
\hspace*{\fill} $\square$(Claim~\ref{nonsatellite})

\medskip

It follows that $K_{\overline{C}}^p$ is a hyperbolic knot in $S^3$.  

Since $K_{\overline{C}}^p$ is a hyperbolic knot with period $p > 2$,  
it follows from \cite[Corollary 1.4]{MM5} that 
$K_{\overline{C}}^p(r)$ is a hyperbolic $3$--manifold for all $r \in \mathbb{Q}$, 
or $p = 3$, $r = 0$ and $g(K_{\overline{C}}^p) =1$. 
Since $g(K_{\overline{C}}^p) \ge p g(\overline{K}) \ge p > 2$, 
the latter cannot occur. 
Hence $K_{\overline{C}}^p(r)$ is a hyperbolic $3$--manifold for all $r \in \mathbb{Q}$ as desired. 

\medskip

$(2)$ 
In the following, 
for notational simplicity, 
we write $K_i = K_{\overline{C_i}}^p$. 

The assumption, together with $(1)$,  
implies that $K_i$ $(i = 1, 2)$ is a hyperbolic knot. 
Recall that $K_i$ has an orientation preserving diffeomorphism $f_i$ such that 
$f_i(K_i) = K_i$, $f_i^p = id$ and $\mathrm{Fix}(f_i) = C_i$. 
Note that 
$\overline{K} = K_i / \langle f_i \rangle$ and 
$\overline{C_i} = C_i / \langle f_i \rangle$. 
Suppose that $K_1$ and $K_2$ are isotopic in $S^3$. 
Then we have an orientation preserving diffeomorphism $\varphi$ of $S^3$ such that 
$\varphi(K_1) = K_2$. 
Note that $f'_2 = \varphi^{-1}\circ f_2 \circ \varphi$ is an orientation preserving diffeomorphism of $S^3$, 
which satisfies $f'_2(K_1) = K_1$ and 
gives also a cyclic period $p$ for $K_1$. 
Let us put $C'_2 = \varphi^{-1}(C_2)$. 
Then we see that $\mathrm{Fix}(f'_2) = C'_2$. 
Since $\varphi \circ f'_2 = f_2 \circ \varphi$, 
$\varphi$ induces an orientation preserving diffeomorphism 
$\overline{\varphi} : S^3 = S^3 / \langle f'_2 \rangle \to S^3 = S^3 / \langle f_2 \rangle$ 
sending $K_1 / \langle f'_2 \rangle$ to $\overline{K} = K_2 / \langle f_2 \rangle$ and 
$C'_2 /  \langle f'_2 \rangle$ to $\overline{C_2} = C_2 /  \langle f_2 \rangle$. 

Now the hyperbolic knot $K_1$ has two orientation preserving, periodic diffeomorphisms $f_1$ and $f'_2$ of period $p > 2$. 
Then \cite[2.1 Theorem (a)]{BZ} shows that 
the pairwise isotopy classes $[f_1]$ and $[f_2' ]$ in the symmetry group $\mathrm{Sym}(S^3, K_1)$ have order $p$. 
Furthermore, since $K_1$ is hyperbolic, 
$\mathrm{Sym}(S^3, K_1)$ is isomorphic to a finite cyclic group or a dihedral group 
\cite[Theorems 10.5.3 and 10.6.2(2)]{Ka}.   
This implies that subgroups $\langle [f_1] \rangle$ and $\langle [f'_2] \rangle$ of order $p$ coincide, 
since $p > 2$.   
Then it follows from \cite[2.1 Theorem (c)]{BZ} that 
$\langle f_1 \rangle$ and $\langle f'_2 \rangle$ are conjugate by a diffeomorphism $g$ in 
$\mathrm{Diff}(S^3, K_1)$ which is isotopic to the identity. 
Hence $(f'_2)^k = g \circ f_1 \circ g^{-1}$ for some integer $k$ $(1 \le k \le p-1)$, 
which has also period $p$. 
Note that $(f'_2)^k$ leaves $K_1$ invariant and $K_1 /  \langle (f'_2)^k \rangle = K_1 / \langle f'_2 \rangle$, 
and that $\mathrm{Fix}((f'_2)^k) = C'_2$ and $C'_2 / \langle (f'_2)^k \rangle = C'_2 / \langle f'_2 \rangle$. 
For any $x \in C_1 = \mathrm{Fix}(f_1)$, 
we have $(f'_2)^k(g(x)) = g(f_1(x)) = g(x)$, 
thus $g(x) \in \mathrm{Fix}((f'_2)^k) = C'_2$, and hence $g(C_1) \subset C'_2$. 
Conversely if $x' \in C'_2 =  \mathrm{Fix}((f'_2)^k)$, 
then we see that $g^{-1}(x') \in C_1$ and $x' \in g(C_1)$, hence $C'_2 \subset g(C_1)$. 
Thus we have $g(C_1) = C'_2$. 
Therefore we have an orientation preserving diffeomorphism 
$\overline{g} : S^3 = S^3 / \langle f_1 \rangle \to S^3 = S^3 /  \langle (f'_2)^k \rangle = S^3 / \langle f'_2 \rangle$  
sending $\overline{K} = K_1 /  \langle f_1 \rangle$ to $K_1 /  \langle (f'_2)^k \rangle = K_1 / \langle f'_2 \rangle$ 
and $\overline{C_1} = C_1 /  \langle f_1 \rangle$ to $C'_2 / \langle (f'_2)^k \rangle = C'_2 / \langle f'_2 \rangle$. 

Now the orientation preserving diffeomorphism 
$\overline{\varphi} \circ \overline{g}$ of $S^3$ satisfies 
$\overline{\varphi} \circ \overline{g}(\overline{K}) = \overline{K}$ and 
$\overline{\varphi} \circ \overline{g}(\overline{C_1}) = \overline{C_2}$. 
Thus $\overline{K} \cup \overline{C_1}$ and $\overline{K} \cup \overline{C_2}$ are isotopic.  
\hspace*{\fill} $\square$(Theorem~\ref{hyperbolic axis2})

\medskip

\section{Examples}
\label{examples}

In this section, 
we present two examples illustrating how the periodic construction works according to whether the initial knot $\overline{K}$ 
is fibered or not fibered. 

\medskip

First we apply Theorem~\ref{construction} in the case where $\overline{K}$ is not fibered. 
In such a case we can choose $\overline{C}$ arbitrarily with $lk(\overline{K}, \overline{C}) \ne 0$ to obtain a knot  
$K_{\overline{C}}^p$ having properties $(1)$ and $(2)$ in Theorem~\ref{construction}. 

Let $T_n$ $(n \ne 0, \pm 1)$ be the twist knot illustrated in Figure~\ref{twistknotTn}. 

\begin{figure}[htbp]
\begin{center}
\includegraphics[width=0.25\linewidth]{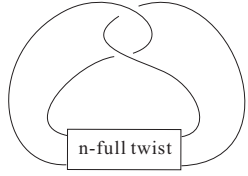}
\caption{The twist knot $T_n$}
\label{twistknotTn}
\end{center}
\end{figure}

Then $T_n$ is a hyperbolic knot, 
and since the Alexander polynomial of $T_n$ is not monic, 
it is not fibered \cite[8.16 Proposition]{BuZ}. 
Suppose that $n > 1$. 
Then it follows from \cite{Tran2, HT2} that 
$\pi_1(T_n(r))$ is left-orderable for $r \in (-4n, 4)$. 
Furthermore, it is known by \cite{Tera} that $\pi_1(T_n(4))$ is left-orderable. 
Hence $\mathcal{S}_{LO}(T_n) \supset  (-4n, 4] \cap \mathbb{Q}$. 

\medskip

\begin{example}
\label{exampleT2_hyperbolic}
Let us take a $2$--component link $T_2 \cup \overline{C}$ as in Figure~\ref{factorknot_nonfibered2}; 
$lk(T_2, \overline{C}) = 1$. 
Let $p$ be any integer with $p > 2$ and $K_{2, \overline{C}}^p$ a knot obtained from $( T_2, \overline{C} )$ by $p$--periodic construction. 
Then $K_{2, \overline{C}}^p$ enjoys the following properties: 
\begin{enumerate}
\item
$K_{2, \overline{C}}^p$ is a hyperbolic knot in $S^3$. 
\item
$K_{2, \overline{C}}^p(r)$ is a hyperbolic $3$--manifold for all $r \in \mathbb{Q}$. 
\item 
$\mathcal{S}_{LO}(K_{2, \overline{C}}^p) \supset (-8p, 4p] \cap \mathbb{Q}$. 
\item 
$\mathcal{S}_{L}(K_{2, \overline{C}}^p) = \emptyset$. 
\end{enumerate}

\begin{figure}[htbp]
\begin{center}
\includegraphics[width=0.3\linewidth]{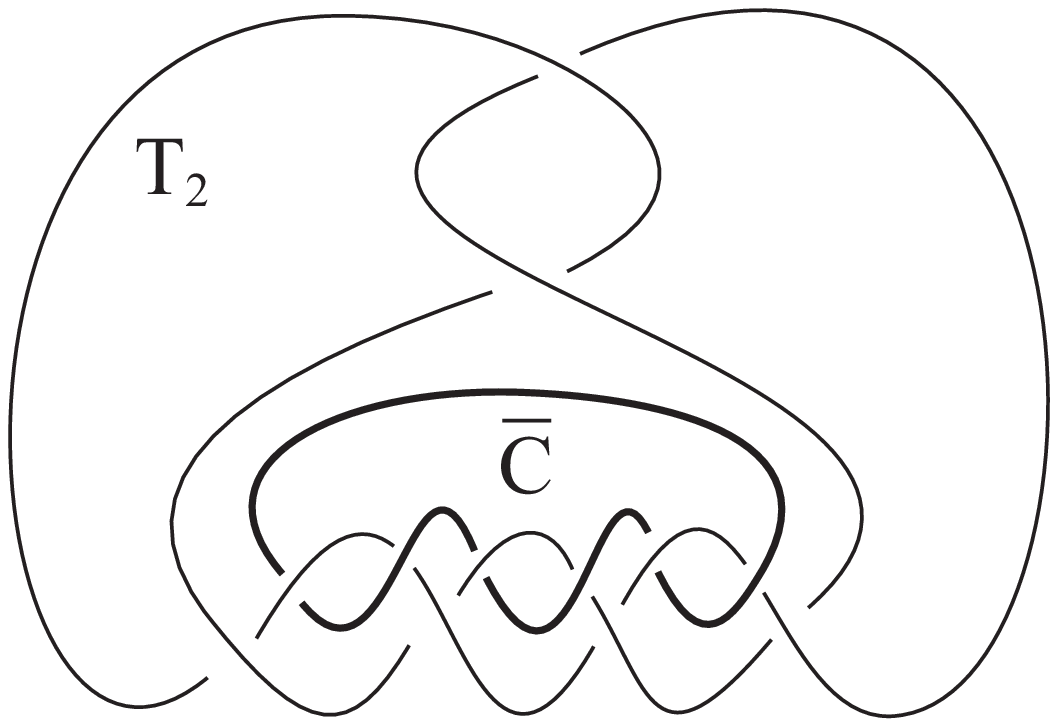}
\caption{The twist knot $T_2$ and an axis $\overline{C}$}
\label{factorknot_nonfibered2}
\end{center}
\end{figure}

\end{example}

\noindent
\textit{Proof.}
Assertions $(1)$ and $(2)$ follow from Theorem~\ref{hyperbolic axis2}(1) 
once we show that $T_2 \cup \overline{C}$ is a hyperbolic link. 
Since $T_2 \cup \overline{C}$ is a non-split prime alternating link \cite[Theorem 1]{Me}, 
it is either a torus link or a hyperbolic link \cite[Corollary 2]{Me}. 
The former cannot happen, because $T_2$ is nontrivial, but $\overline{C}$ is trivial.  
Hence $T_2 \cup \overline{C}$ is a hyperbolic link as desired. 
Since $T_2$ is not fibered and $\pi_1(T_2(r))$ is left-orderable for $r \in (-8, 4]$, 
assertions $(3)$ and $(4)$ follow from Theorem~\ref{construction}.  
\hspace*{\fill} $\square$(Example~\ref{exampleT2_hyperbolic})

\bigskip

Next we apply Theorem~\ref{construction} in the case where $\overline{K}$ is a fibered knot. 
In the next example 
we take a trefoil knot $T_{-3, 2}$ as $\overline{K}$. 

\medskip

\begin{example}
\label{example trefoil}
Let us take the $2$--component link $T_{-3, 2} \cup \overline{C}$ shown in Figure~\ref{geom_alg1}; 
$lk(T_{-3, 2}, \overline{C}) = 1$. 
Let $p$ be any integer with $p > 2$ and 
$K_{-3, 2, \overline{C}}^p$ a knot obtained from $( T_{-3, 2}, \overline{C} )$ by $p$--periodic construction. 
Then $K_{-3, 2, \overline{C}}^p$ enjoys the following properties: 
\begin{enumerate}
\item
$K_{-3, 2, \overline{C}}^p$ is a hyperbolic knot in $S^3$. 
\item
$K_{-3, 2, \overline{C}}^p(r)$ is a hyperbolic $3$--manifold for all $r \in \mathbb{Q}$.
\item 
$\mathcal{S}_{LO}(K_{-3, 2, \overline{C}}^p) \supset (-p, \infty) \cap \mathbb{Q}$. 
\item 
$\mathcal{S}_{L}(K_{-3, 2, \overline{C}}^p) = \emptyset$. 
\end{enumerate}

\begin{figure}[htbp]
\begin{center}
\includegraphics[width=0.7\linewidth]{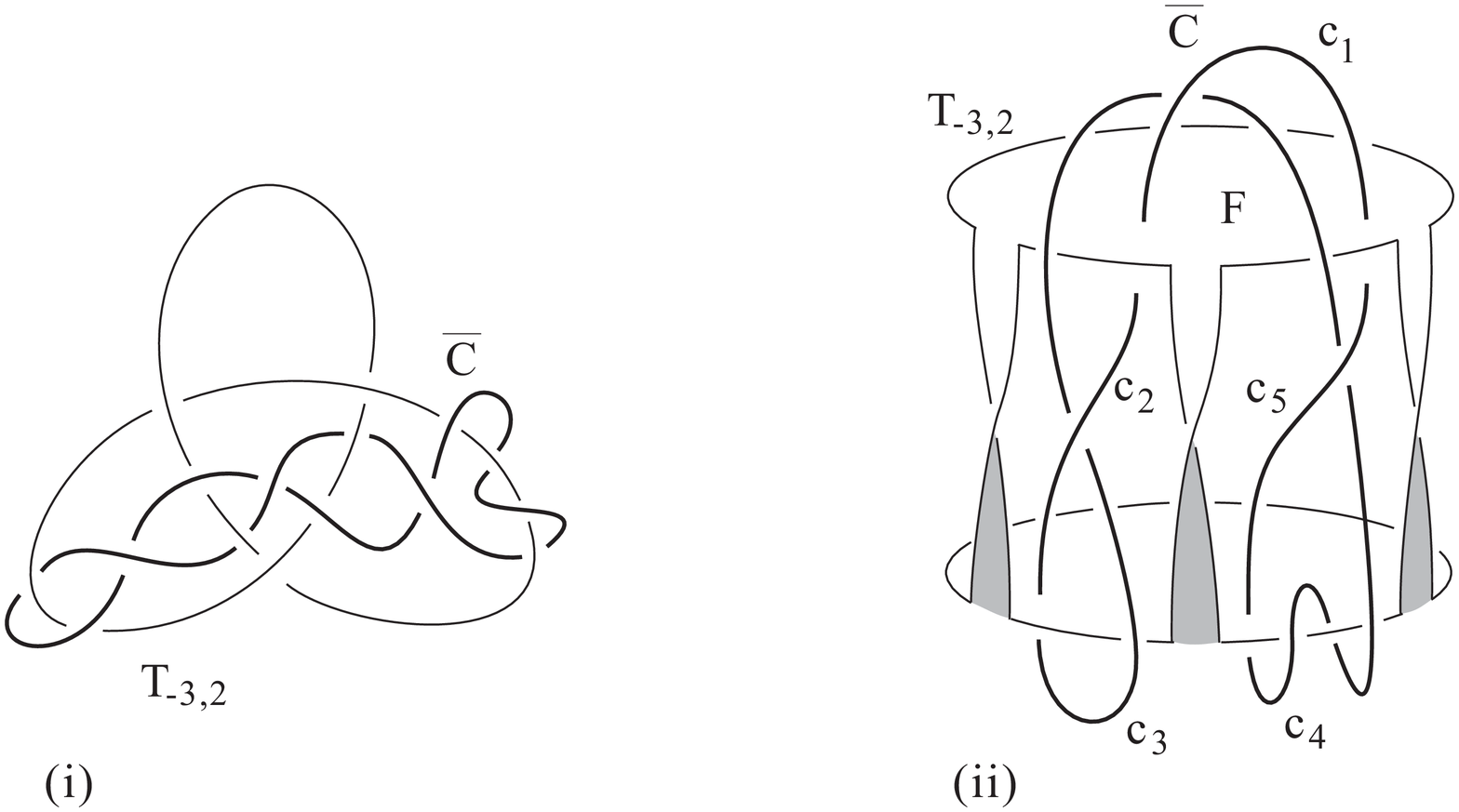}
\caption{The trefoil knot $T_{-3, 2}$ and the unknotted circle $\overline{C}$}
\label{geom_alg1}
\end{center}
\end{figure}

\end{example}

\noindent
\textit{Proof of Example~\ref{example trefoil}.}
Recall that  $\mathcal{S}_{LO}(T_{-3, 2}) = (-1, \infty) \cap \mathbb{Q}$; 
see Remark~\ref{mirror}(2) and Example~\ref{torus knots}. 

Since as illustrated in Figure~\ref{geom_alg1}(i) 
$T_{-3, 2} \cup \overline{C}$ is a non-split prime alternating link \cite[Theorem 1]{Me}, 
it is either a torus link or a hyperbolic link \cite[Corollary 2]{Me}. 
If we have the former case, 
then $T_{-3, 2}$ is isotopic to $\overline{C}$ which is a trivial knot, a contradiction. 
Hence $T_{-3, 2} \cup \overline{C}$ is a hyperbolic link. 
Then $(1)$ and $(2)$ follow from Theorem~\ref{hyperbolic axis2}(1). 

Let us prove $(3)$ and $(4)$ using Theorem~\ref{construction}. 
Since $T_{-3, 2}$ is fibered, 
we need to check the condition of Theorem~\ref{construction}: 
for any fiber surface $\overline{S}$ of $E(T_{-3, 2})$, 
$|\overline{S} \cap \overline{C} |$ is strictly bigger than 
the algebraic intersection number between $\overline{S}$ and $\overline{C}$, 
i.e. $lk(T_{-3, 2}, \overline{C})$. 

In Figure~\ref{geom_alg1}(ii), 
we give a minimal genus Seifert surface $F$ of $T_{-3, 2}$, 
which is a once-punctured torus with $\partial F = T_{-3, 2}$.  
Put $\overline{S} = F \cap E(T_{-3, 2})$. 
Then by \cite[Lemma 5.1]{EL} $\overline{S}$ is a fiber surface of $E(T_{-3, 2})$. 
We see that $|\overline{S} \cap \overline{C} | = 5$ and 
the algebraic intersection number between $\overline{S}$ and $\overline{C}$ is one. 
Assume for a contradiction that we have another fiber surface $\overline{S}'$ of $E(T_{-3, 2})$ 
such that $|\overline{S}' \cap \overline{C} | < |\overline{S} \cap \overline{C} |$. 
Since $\overline{S}$ and $\overline{S}'$ are fiber surfaces of $E(T_{-3, 2})$, 
they are isotopic; see \cite[Lemma 5.1]{EL}, \cite{Thu}. 
This then implies that 
we can isotope $\overline{C}$ to $\overline{C}'$ in $E(T_{-3, 2})$ so that 
$|\overline{S} \cap \overline{C}'  | < |\overline{S} \cap \overline{C} |$. 

\begin{claim}
\label{semidisk}
There exists a smooth map $\varphi$ from a semi-disk $D$ into $E(T_{-3, 2})$ 
such that $\varphi^{-1}(\overline{C})$ is an arc $c \subset \partial D$ and 
$\varphi^{-1}(\overline{S})$ is the arc $\alpha =\overline{\partial D - c}$.  
\end{claim}

\noindent
\textit{Proof of Claim~\ref{semidisk}.}
Let $\Phi : S^1 \times [0, 1] \to E(T_{-3, 2})$ be a smooth map giving an isotopy between 
$\overline{C} (= \Phi(S^1 \times \{ 0 \}))$ to $\overline{C}' (= \Phi(S^1 \times \{ 1 \}))$. 
We may assume $\Phi$ is transverse to $\overline{S}$. 
Furthermore, the essentiality of $\overline{S}$ in $E(T_{-3, 2})$ enables us to modify  $\Phi$ to 
eliminate the circle components as usual.  
Since $|\overline{S} \cap \overline{C}'  | < |\overline{S} \cap \overline{C} | = 5$ 
and the algebraic intersection number between 
$\overline{S}$ and $\overline{C}'$ coincides with the algebraic intersection number between 
$\overline{S}$ and $\overline{C}$,  
we have $|\overline{S} \cap \overline{C}'  | = 1$ or $3$. 
Thus $\Phi^{-1}(\overline{S})$ consists of 
three properly embedded arcs $\alpha,\  \alpha'$ and $\beta$, 
where $\partial \alpha \subset S^1 \times \{ 0 \}$, 
$\partial \alpha' \subset S^1 \times \{ 0 \}$, 
and $\beta$ connects $S^1 \times \{ 0 \}$ and $S^1 \times \{ 1 \}$ (Figure~\ref{annulus}(i), (ii)),  
consists of four properly embedded arcs $\alpha,\ \beta,\ \beta'$ and $\beta''$, 
where $\partial \alpha \subset S^1 \times \{ 0 \}$, 
and each of $\beta, \beta', \beta''$ connects $S^1 \times \{ 0 \}$ and $S^1 \times \{ 1 \}$ (Figure~\ref{annulus}(iii)), 
or consists of four properly embedded arcs $\alpha,\ \alpha', \beta$ and $\gamma$, 
where $\partial \alpha \subset S^1 \times \{ 0 \}$, 
$\partial \alpha' \subset S^1 \times \{ 0 \}$, 
$\beta$ connects $S^1 \times \{ 0 \}$ and $S^1 \times \{ 1 \}$, 
and $\partial \gamma \subset S^1 \times \{ 1 \}$  (Figure~\ref{annulus}(iv), (v)). 
In either case there is a semi-disk $D$ cobounded by $\alpha$ and an arc $c \subset S^1 \times \{ 0 \}$.  

\begin{figure}[htbp]
\begin{center}
\includegraphics[width=1.0\linewidth]{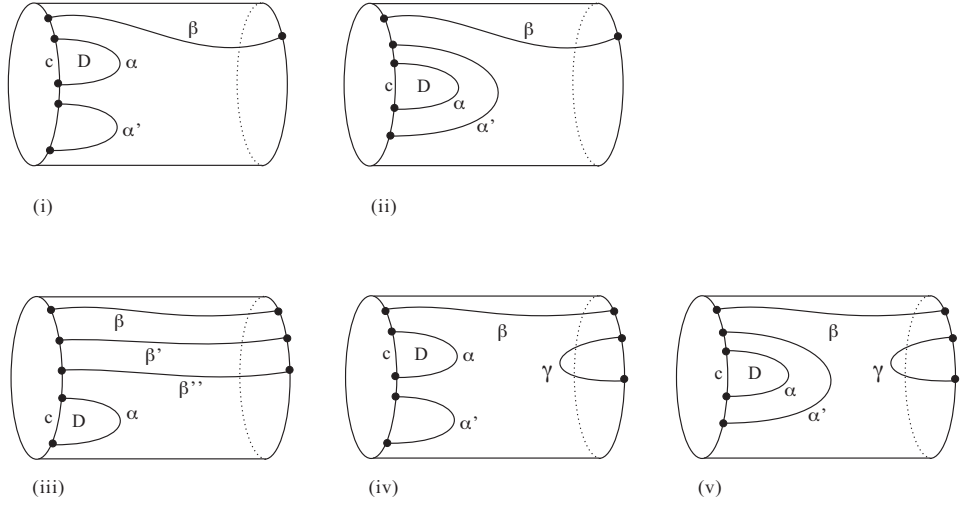}
\caption{$\Phi^{-1}(\overline{S})$ in $S^1 \times [0, 1]$}
\label{annulus}
\end{center}
\end{figure}

Putting $\varphi = \Phi |_D : D \to E(T_{-3, 2})$, 
we obtain a desired smooth map. 
\hspace*{\fill} $\square$(Claim~\ref{semidisk})

\medskip

Cut open $E(T_{-3, 2})$ along $\overline{S}$ to obtain a product $3$--manifold $\overline{S} \times [0, 1]$. 
The circle $\overline{C}$ is cut into five arcs $c_1, c_2, c_3, c_4$ and $c_5$ as in Figure~\ref{geom_alg1}(ii).  
Note that $\partial c_1 \subset \overline{S} \times \{ 0 \}$, $\partial c_3 \subset \overline{S} \times \{ 1 \}$, 
and each of $c_2, c_4, c_5$ connects $\overline{S} \times \{ 0 \}$ and $\overline{S} \times \{ 1 \}$. 
Moreover, 
we see that  $c_1$ and $c_3$ are linking once relative their boundaries. 

On the other hand, 
since $c$ is either $c_1$ or $c_3$, 
Claim~\ref{semidisk} shows that $c_1$ and $c_3$ are unlinked relative their boundaries. 
This contradiction shows that for any fiber surface $\overline{S}$, 
$| \overline{S} \cap \overline{C}  | = 5$ and 
$|  \overline{S} \cap \overline{C} | > lk(T_{-3, 2}, \overline{C})$.  

Since $\pi_1(T_{-3, 2}(r))$ is left-orderable if $r \in (-1, \infty)$, 
the conclusions $(3)$ and $(4)$ follow from Theorem~\ref{construction}. 
This completes the proof of Example~\ref{example trefoil}. 
\hspace*{\fill} $\square$(Example~\ref{example trefoil})

\bigskip

\section{Surgeries on alternating knots}
\label{periodic_alternating}

Theorem 1.5 in \cite{OS3}, 
together with \cite[Proposition 9.6]{OS4} (\cite[Proof of Corollary 1.3]{Ni}, \cite[Claim 2]{IJ}), 
shows that for an alternating knot $K$ which is not a $(p, 2)$--torus knot, 
$K(r)$ is not an $L$--space for all $r \in \mathbb{Q}$. 

We say that an alternating knot is \textit{positive} (resp. \textit{negative}) if 
it has a reduced alternating diagram such that 
each of the crossings is positive (resp. negative). 
An alternating knot is \textit{special} if it is either positive or negative. 

\medskip

In \cite{BGW} 
Boyer, Gordon and Watson prove:  

\begin{proposition}[\cite{BGW}]
\label{BGWalternating}
Let $K$ be a prime alternating knot in $S^3$. 
 
\begin{enumerate}
\item
If $K$ is not a special alternating knot, 
then $\pi_1(K(\frac{1}{n}))$ is left-orderable for all non-zero integers $n$. 
\item
If $K$ is a positive $($resp. negative$)$ alternating knot, 
then $\pi_1(K(\frac{1}{n}))$ is left-orderable for all positive $($resp. negative$)$ integers $n$.  
\end{enumerate}
\end{proposition}

Let $\overline{K}$ be an alternating knot. 
For convenience, 
we position $\overline{K} \subset \mathbb{R}^3 = S^3 - \{ \infty \}$ so that 
$\overline{K}$ lies in the $xy$--plane except near crossings of $\overline{K}$, 
where $\overline{K}$ lies on a ``bubble" as in \cite{Me}. 
Then we say an unknotted circle $\overline{C} \subset S^3 - \overline{K}$ is \textit{perpendicular} if it passes $\infty$ and intersects 
the $xy$--plane exactly once. 
Note that $\overline{C} \cap \mathbb{R}^3$ is perpendicular to the $xy$--plane.  
See Figure~\ref{alternating}, 
in which the dot indicates a perpendicular circle $\overline{C}$. 

\begin{figure}[htbp]
\begin{center}
\includegraphics[width=0.23\linewidth]{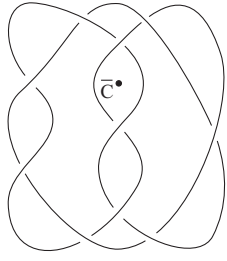}
\caption{An alternating knot $\overline{K}$ and a perpendicular circle $\overline{C}$}
\label{alternating}
\end{center}
\end{figure}

\begin{proposition}
\label{p_alternating}
Let $\overline{K}$ be a prime alternating knot and $\overline{C}$ a perpendicular circle. 
Let $p$ be an integer such that $p \ge 2$ and $(p, lk(\overline{K}, \overline{C})) = 1$,  
and let $K_{\overline{C}}^p$ be the knot obtained from $( \overline{K}, \overline{C} )$ by $p$--periodic construction. 
Then we have: 
 
\begin{enumerate}
\item
$K_{\overline{C}}^p$ is an alternating knot. 
\item
If $\overline{K}$ is not a special alternating knot, 
then $\pi_1(K_{\overline{C}}^p(\frac{p}{n}))$ is left-orderable for all non-zero integers $n$. 
\item
If $\overline{K}$ is a positive $($resp. negative$)$ alternating knot, 
then 
$\pi_1(K_{\overline{C}}^p(\frac{p}{n}))$ is left-orderable for all positive $($resp. negative$)$ integers $n$. 
\end{enumerate}
\end{proposition}

\medskip

\noindent
\textit{Proof of Proposition~\ref{p_alternating}.}
The first assertion follows immediately from diagramatic consideration. 
The conclusions $(2)$ and $(3)$ follow from 
Proposition~\ref{BGWalternating} and Theorem~\ref{factorknot_LO}. 
\hspace*{\fill} $\square$(Proposition~\ref{p_alternating})

\medskip

\begin{remark}
\label{nonfibered_alternating}
In Proposition~\ref{p_alternating}, 
if $\overline{K}$ is not a $(p, 2)$--torus knot, 
then $K_{\overline{C}}^p$ is not a $(p', 2)$--torus knot. 
For otherwise, 
the argument in the proof of Lemma~\ref{irreducible} implies that 
$\overline{K}$ is a torus knot. 
Since it is alternating, it is a $(p, 2)$--torus knot for some odd integer $p$ \cite[Theorem 3.2]{Mu}, 
a contradiction. 
Therefore, 
as mentioned in the beginning of this section, 
$K_{\overline{C}}^p(r)$ is not an $L$--space for all $r \in \mathbb{Q}$. 
\end{remark}

\medskip 

Applying Proposition~\ref{p_alternating} and Remark~\ref{nonfibered_alternating}, 
we have: 

\medskip

\begin{example}
\label{10_61}
Take an alternating knot $\overline{K}$  
and a perpendicular circle $\overline{C}$ as illustrated in Figure~\ref{alternating}; $lk(\overline{K}, \overline{C}) = 1$. 
Note that $\overline{K}$ is not a special alternating knot. 
Hence for any integer $p \ge 2$, 
$K_{\overline{C}}^p$ is an alternating knot,  
$K_{\overline{C}}^p(r)$ is not an $L$--spaces for all $r \in \mathbb{Q}$, 
and 
$\pi_1(K_{\overline{C}}^p(\frac{p}{n}))$ is left-orderable for all non-zero integers $n$. 
\end{example}

\bigskip

\section{Knots with $\mathcal{S}_{LO}(K) = \mathbb{Q}$ and $\mathcal{S}_{L}(K) = \emptyset$}
\label{proofs}

The goal of this section is to prove Theorems~\ref{satellite} and \ref{all rationals}. 
We start with Proposition~\ref{ClayWatson} below, 
which was shown by Clay and Watson \cite[Proposition~4.1]{CW1}. 

Let $k$ be a knot in $S^3$, 
which is contained in a standardly embedded solid torus $V \subset S^3$. 
Assume that $k$ is not contained in a $3$--ball in $V$. 
We call $k$ a \textit{pattern knot} in $S^3$ and the pair $(V, k)$ a \textit{pattern}. 
Let $f$ be an orientation preserving embedding from $V$ into $S^3$ which sends a preferred longitude of $V$  to that of 
$f(V) \subset S^3$. 
Then we obtain a knot $K = f(k)$ in $S^3$, which is called a \textit{satellite knot} with a pattern knot $k$ and 
a \textit{companion knot} $K' = f(c)$, where $c$ is a core of $V$. 

\begin{proposition}[\cite{CW1}]
\label{ClayWatson}
Let $K$ be a satellite knot with a pattern knot $k$. 
If $K(r)$ is irreducible and $r \in \mathcal{S}_{LO}(k)$, 
then $r \in \mathcal{S}_{LO}(K)$. 
\end{proposition}

\subsection{Composite knots $K$ with $\mathcal{S}_{LO}(K) = \mathbb{Q}$ and $\mathcal{S}_{L}(K) = \emptyset$}

In this subsection 
we prove that the connected sum of two torus knots $T_{-p, q}$ and $T_{r, s}$ where $p > q \ge 2$ and $r > s \ge 2$, 
satisfies $\mathcal{S}_{LO}(T_{-p, q}\, \sharp\, T_{r, s}) = \mathbb{Q}$ and 
$\mathcal{S}_{L}(T_{-p, q}\, \sharp\, T_{r, s}) = \emptyset$ (Proposition~\ref{T-pqTrs}).  
Thus $T_{-p, q}\, \sharp\, T_{r, s}$ satisfies Conjecture~\ref{LspaceConjecture2}.

Proposition~\ref{ClayWatson} and Theorem~\ref{factorknot_LO} immediately imply: 

\begin{proposition}
\label{composite}
Let $K$ and $K'$ be nontrivial knots. 
Then we have: 
\begin{enumerate}
\item
$\mathcal{S}_{LO}(K\, \sharp\, K') \supset \mathcal{S}_{LO}(K) \cup \mathcal{S}_{LO}(K')$. 

\item
$\mathcal{S}_{LO}(pK) \supset p \mathcal{S}_{LO}(K)$, where $pK$ denotes the connected sum of $p$ copies of $K$. 
\end{enumerate}
\end{proposition}

\noindent
\textit{Proof of Proposition~\ref{composite}.}
$(1)$ Following \cite[Lemma~7.1]{Go}, 
we see that $(K\, \sharp\, K')(r)$ is irreducible for all $r \in \mathbb{Q}$. 

Let us regard $K\, \sharp\, K'$ as a satellite knots with a pattern knot $K$ and a companion knot $K'$. 
Then Proposition~\ref{ClayWatson} 
shows that $\mathcal{S}_{LO}(K\, \sharp\, K') \supset \mathcal{S}_{LO}(K)$. 
Exchanging the roles of $K$ and $K'$,  
we have $\mathcal{S}_{LO}(K\, \sharp\, K') \supset \mathcal{S}_{LO}(K')$ as well. 
Thus $\mathcal{S}_{LO}(K\, \sharp\, K') \supset \mathcal{S}_{LO}(K) \cup \mathcal{S}_{LO}(K')$. 

\medskip 

$(2)$   
Since $pK$ is a knot with cyclic period $p$ whose factor knot is $K$, 
the result follows from Theorem~\ref{factorknot_LO}. 
\hspace*{\fill} $\square$(Proposition~\ref{composite})

\medskip

As a step toward proofs of Theorems~\ref{satellite} and \ref{all rationals}, 
we prove: 

\begin{proposition}
\label{T-pqTrs}
For torus knots $T_{-p, q}$ and $T_{r, s}$, 
where $p > q \ge 2$ and $r > s \ge 2$, 
we have $\mathcal{S}_{LO}(T_{-p, q}\, \sharp\, T_{r, s}) = \mathbb{Q}$ and 
$\mathcal{S}_{L}(T_{-p, q}\, \sharp\, T_{r, s}) = \emptyset$. 
\end{proposition}

\noindent
\textit{Proof of Proposition~\ref{T-pqTrs}.}
Recall that $\mathcal{S}_{LO}(T_{-p, q}) = (-pq+p+q, \infty) \cap \mathbb{Q}$ and 
$\mathcal{S}_{LO}(T_{r, s}) = (-\infty, rs-r-s) \cap \mathbb{Q}$. 
Note that $-pq+p+q < 0 < rs-r-s$. 
Now apply Proposition~\ref{composite} to $T_{-p, q}\, \sharp\, T_{r, s}$ to conclude that 
$\mathcal{S}_{LO}(T_{-p, q}\, \sharp\, T_{r, s}) \supset ((-pq+p+q, \infty) \cup (-\infty, rs-r-s)) \cap \mathbb{Q} = \mathbb{Q}$. 
Hence $\mathcal{S}_{LO}(T_{-p, q}\, \sharp\, T_{r, s}) = \mathbb{Q}$.  

Next we show that $T_{-p, q}\, \sharp\, T_{r, s}$ has no $L$--space surgeries.  

\begin{claim}
\label{Alexander}
The coefficient of $t$ in the Alexander polynomial of 
$T_{-p, q}\, \sharp\, T_{r, s}$ is $-2$. 
\end{claim}

\noindent
\textit{Proof of Claim~\ref{Alexander}.}
Recall that $T_{-p, q}$ has the Alexander polynomial 
$\displaystyle \Delta_{T_{-p, q}}(t) = \Delta_{T_{p, q}}(t) = \frac{(t^{pq} - 1)(t-1)}{(t^p -1)(t^q -1)}$, 
and $T_{r, s}$ has the Alexander polynomial 
$\displaystyle \Delta_{T_{r, s}}(t) = \frac{(t^{rs} - 1)(t-1)}{(t^r -1)(t^s -1)}$. 
Since $\Delta_{T_{-p, q}\, \sharp\, T_{r, s}}(t) = \Delta_{T_{-p, q}}(t) \Delta_{T_{r, s}}(t)
= \Delta_{T_{p, q}}(t) \Delta_{T_{r, s}}(t)$, \break
$\Delta'_{T_{-p, q}\, \sharp\, T_{r, s}}(0) = 
\Delta'_{T_{p, q}}(0) \Delta_{T_{r, s}}(0) +  
\Delta_{T_{p, q}}(0)  \Delta'_{T_{r, s}}(0)$. 
Note that $\Delta_{T_{p, q}}(0) = \Delta_{T_{r, s}}(0) = 1$ and 
a simple computation shows that 
$ \Delta'_{T_{p, q}}(0) =  \Delta'_{T_{r, s}}(0) = -1$. 
Thus $\Delta'_{T_{-p, q}\, \sharp\, T_{r, s}}(0) =  (-1) + (-1) = -2$. 
This then implies that 
the coefficient of $t$ in the Alexander polynomial of 
$T_{-p, q}\, \sharp\, T_{r, s}$ is $-2$. 
\hspace*{\fill} $\square$(Claim~\ref{Alexander})

\medskip

Apply \cite[Corollary~1.3]{OS3}, together with \cite[Proposition 9.6]{OS4} (\cite[Proof of Corollary 1.3]{Ni}, \cite[Claim 2]{IJ}), 
to conclude that $T_{-p, q}\, \sharp\, T_{r, s}$ has no $L$--space surgeries. 

\hspace*{\fill} $\square$(Proposition~\ref{T-pqTrs})

\medskip 

Let us consider 
the connected sum $T_{p, q}\, \sharp\, T_{r, s}$ 
instead of $T_{-p, q}\, \sharp\, T_{r, s}$, 
where $p > q \ge 2$ and $r > s \ge 2$. 
The argument in the proof of Claim~\ref{Alexander} shows that $\mathcal{S}_{L}(T_{p, q}\, \sharp\, T_{r, s}) = \emptyset$. 
On the other hand, 
putting $m_0 = \mathrm{max} \{ pq -p-q,\ rs - r -s \}$, 
Example~\ref{torus knots} and Proposition~\ref{composite} merely imply 
$\mathcal{S}_{LO}(T_{p, q}\, \sharp\, T_{r, s}) \supset (-\infty, m_0)$. 
So we would like to ask: 

\begin{question}
\label{T-pqTrs_LO}
Does $\mathcal{S}_{LO}(T_{p, q}\, \sharp\, T_{r, s}) = \mathbb{Q}$ hold for integers $p > q \ge 2$ and $r > s \ge 2$? 
\end{question}

\smallskip

\subsection{Proof of Theorem~\ref{satellite}}
Let us consider $k = T_{-3, 2} \sharp\, T_{3, 2}$ and 
take an unknotted circle $C$ as in Figure~\ref{2trefoil}. 
Following Proposition~\ref{T-pqTrs}, 
$\mathcal{S}_{LO}(k) = \mathbb{Q}$. 

\begin{figure}[htbp]
\begin{center}
\includegraphics[width=0.6\linewidth]{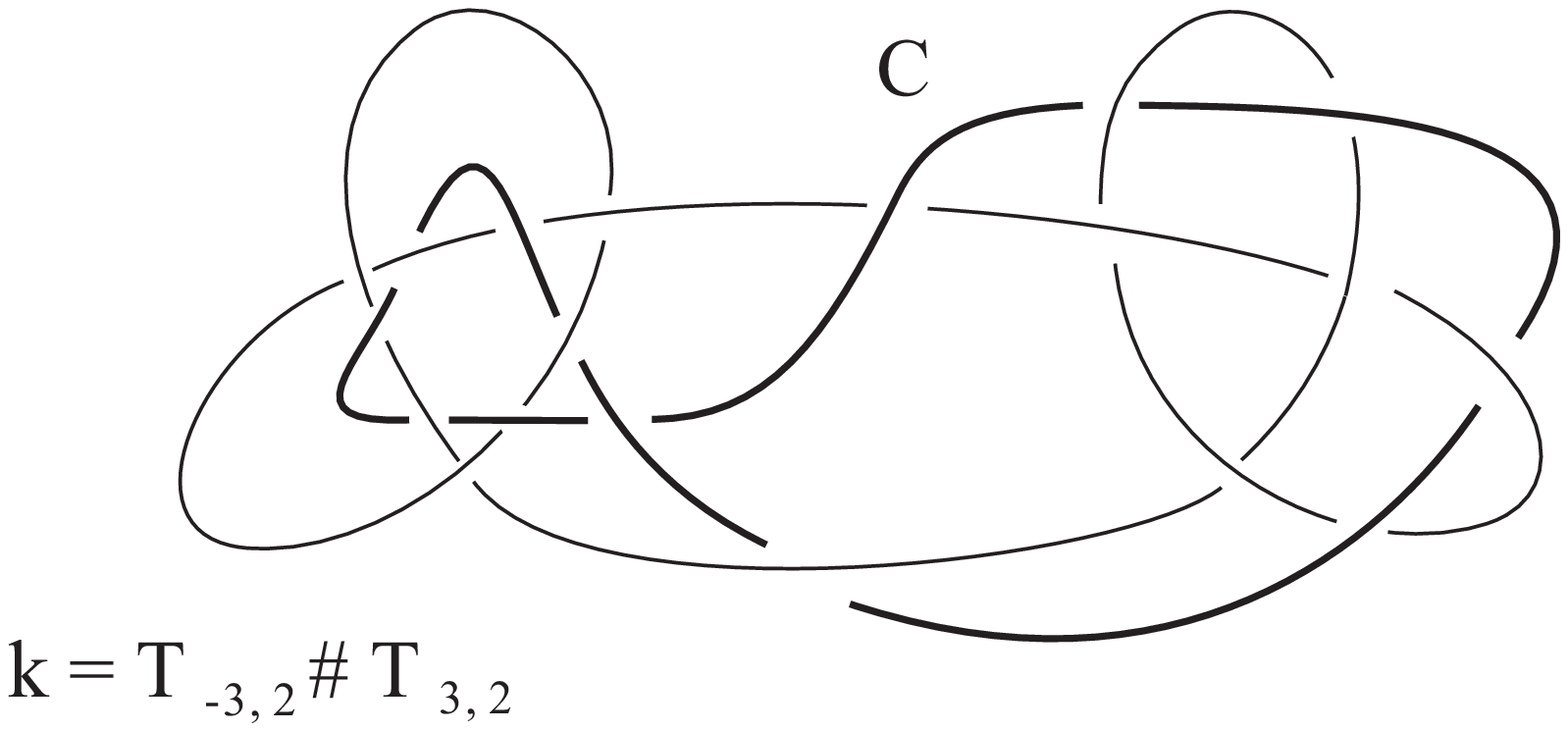}
\caption{$k \cup C$}
\label{2trefoil}
\end{center}
\end{figure}

Note that the link $k \cup C$ is an alternating link. 
Since $k \cup C$ is a non-split prime alternating link 
\cite[Theorem~1]{Me}, 
it is either a torus link or a hyperbolic link \cite[Corollary~2]{Me}. 
The former is not the case, 
because $k$ is nontrivial, but $C$ is trivial. 
Thus $k\cup C$ is hyperbolic, 
hence letting $V = S^3 - \mathrm{int}N(C)$, 
$k$ is a hyperbolic knot in $V$. Apply the satellite construction with 
the pattern $(V, k)$ and the companion knot $K'$ to obtain a satellite knot $K$ 
with a pattern knot $k = T_{-3, 2} \sharp\, T_{3, 2}$. 
Since $k$ is hyperbolic in $V$, 
the satellite knot $K$ is prime, 
and the $3$--manifold obtained from $V$ by $r$--surgery on $k$ 
is again hyperbolic for all but finitely many $r \in \mathbb{Q}$. 
This then implies that $K(r)$ is a toroidal $3$--manifold with a hyperbolic piece,  
in particular, 
$K(r)$ is not a graph manifold, for all but finitely many $r \in \mathbb{Q}$. 
This establishes $(1)$. 

If $K(r)$ were reducible for some $r \in \mathbb{Q}$, 
then $K$ is cabled \cite[4.5 Corollary]{Sch}. 
However, this is impossible, 
because $V- k$ is hyperbolic.  
Hence $K(r)$ is irreducible for any $r \in \mathbb{Q}$. 
Now Proposition~\ref{ClayWatson} shows that 
$\mathcal{S}_{LO}(K) \supset \mathcal{S}_{LO}(k) = \mathbb{Q}$. 

Let us see that $\mathcal{S}_{L}(K) = \emptyset$. 
Since $lk(k, C) = 0$, 
i.e. the winding number of $K$ in $V$ is zero, 
$(V, k)$ is not fibered, 
and hence neither is the satellite knot $K$; 
see \cite[Theorem~1]{HMS}. 
Hence \cite[Corollary~1.3]{Ni} shows that 
$\mathcal{S}_{L}(K) = \emptyset$. 

Finally we show that there are infinitely many satellite knots $K$ with a companion knot $K'$ and enjoy the required properties in Theorem~\ref{satellite}. 
For instance, let us take $k_p = T_{-3, 2} \sharp T_{p, 2}$ $(p \ge 3)$. 
As shown in Proposition~\ref{T-pqTrs}, 
$\mathcal{S}_{LO}(k_p) = \mathbb{Q}$ for all $p \ge 3$. 
It follows from Theorem~\ref{hyperbolic axis1}, 
there is an unknotted circle $C_p$ so that 
$k_p \cup C_p$ is hyperbolic and $lk(k_p,C_p) = 0$. 
Each $C_p$ gives a pattern $(V, k_p)$. 
Let $K_p$ be a satellite knot with a companion knot $K'$ and pattern $(V, k_p)$. 
Then the same argument as above shows that 
$K_p$ satisfies the properties of Theorem~\ref{satellite}. 
If $p \ne p' \ge 3$, then $k_p \cup C_p$ are not isotopic to $k_{p' }\cup C_{p'}$. 
Hence there is no orientation preserving diffeomorphism  of $V$ which leaves the preferred longitude of $V$ invariant and maps 
$k_p$ to $k_{p'}$. 
Thus we see that the resulting satellite knots $K_p$ and $K_{p'}$ are never isotopic. 
\hspace*{\fill} $\square$(Theorem~\ref{satellite})

\medskip

\subsection{Proof of Theorems~\ref{all rationals}}

As in the proof of Theorem~\ref{satellite},  
we consider the connected sum $T_{-3, 2}\, \sharp\, T_{3, 2}$, 
which has the property: 
$\mathcal{S}_{LO}(T_{-3, 2}\, \sharp\, T_{3, 2}) = \mathbb{Q}$ (Proposition~\ref{T-pqTrs}). 

Although we can apply the periodic construction and Theorem~\ref{construction} to the fibered knot 
$T_{-3, 2}\, \sharp\, T_{3, 2}$, 
for ease of handling, 
we take the connected sum $(T_{-3, 2}\, \sharp\, T_{3, 2})\, \sharp\, T_2$, 
where $T_2$ is the twist knot shown in Figure~\ref{twistknotTn}.  
The Alexander polynomial of $(T_{-3, 2}\, \sharp\, T_{3, 2})\, \sharp\, T_2$ is 
$(t^2-t+1)^2(2t^2 - 5t +2)$, which is not monic, 
and hence $(T_{-3, 2}\, \sharp\, T_{3, 2})\, \sharp\, T_2$ is not fibered. 
Proposition~\ref{composite} shows that $\mathcal{S}_{LO}((T_{-3, 2}\, \sharp\, T_{3, 2})\, \sharp\, T_2) 
\supset \mathcal{S}_{LO}(T_{-3, 2}\, \sharp\, T_{3, 2}) = \mathbb{Q}$, 
and hence 
$\mathcal{S}_{LO}((T_{-3, 2}\, \sharp\, T_{3, 2})\, \sharp\, T_2) = \mathbb{Q}$. 

Let us put $\overline{K} = T_{-3, 2}\, \sharp\, T_{3, 2}\, \sharp\, T_2$ and 
take an unknotted circle $\overline{C}$ as in Figure~\ref{2trefoil_twist2}; 
$lk(\overline{K}, \overline{C}) = 1$. 

\begin{figure}[htbp]
\begin{center}
\includegraphics[width=0.85\linewidth]{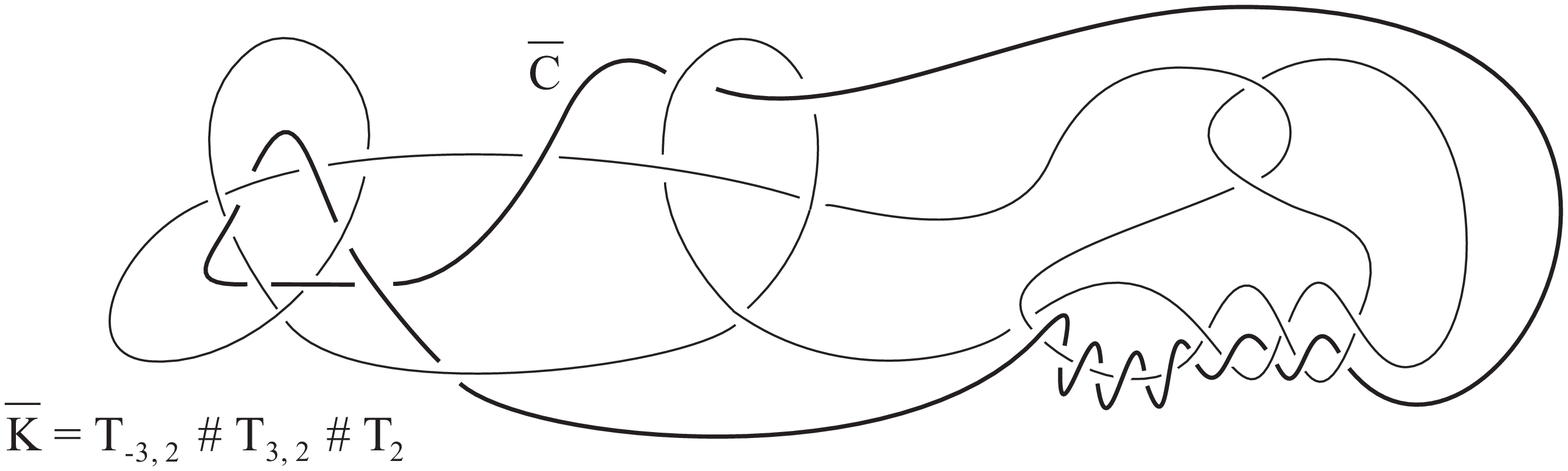}
\caption{$\overline{K} \cup \overline{C}$}
\label{2trefoil_twist2}
\end{center}
\end{figure}

Since $\overline{K} \cup \overline{C}$ is a non-split prime alternating link \cite[Theorem 1]{Me}, 
it is either a torus link or a hyperbolic link \cite[Corollary 2]{Me}. 
The former cannot happen,  
because $\overline{K}$ is nontrivial, but $\overline{C}$ is trivial. 
Hence $\overline{K} \cup \overline{C}$ is a hyperbolic link. 
Let $p > 2$ be any integer, 
and apply the $p$--periodic construction to the pair $(\overline{K}, \overline{C})$ to obtain 
a knot $K_{\overline{C}}^p$. 
It follows from Theorem~\ref{construction} and Theorem~\ref{hyperbolic axis2}$(1)$ that 
$K_{\overline{C}}^p$ is a hyperbolic knot and enjoys the properties $(1), (2)$ and $(3)$ in Theorem~\ref{all rationals}. 
By changing $p$, we obtain infinitely many such knots. 
For instance, see Remark~\ref{genus}. 
\hspace*{\fill} $\square$(Theorem~\ref{all rationals})

\bigskip 

\begin{remark}
\label{infinitely many all rationals}
\begin{enumerate}
\item
By Theorem~\ref{hyperbolic axis1} 
there are infinitely many unknotted circles for $\overline{K} = T_{-3, 2}\, \sharp\, T_{3, 2}\, \sharp\, T_2$, 
and for each unknotted circle $\overline{C}$ we obtain infinitely many hyperbolic knots $K_{\overline{C}}^p$, 
where $p$ and $lk(\overline{K}, \overline{C})$ are relatively prime.  
See also Theorem~\ref{hyperbolic axis2}$(2)$

\item
Recall that any knot $K$ obtained by the periodic construction $($Section~\ref{constructions}$)$,  
for instance a knot obtained in the proof of Theorem~\ref{all rationals},  
is not fibered and every nontrivial surgery on $K$ is a left-orderable, non--$L$--space surgery. 
So we can apply Theorem~\ref{construction} again to the knot $K$ and an arbitrarily chosen unknotted circle to obtain yet further infinitely many non-fibered knots $K'$ each of which 
has the (same) factor knot $K$. 
Then $r$--surgery on $K'$ is also a left-orderable, non--$L$--space surgery for all $r \in \mathbb{Q}$. 
We can apply this procedure repeatedly arbitrarily many times. 

\item
Let $K$ be the knot $10_{99}$ in Rolfsen's knot table \cite{Ro}. 
Recently Clay \cite{C} used an epimorphism from $E(K)$ to $E(T_{3, 2})$ which preserves the peripheral subgroup \cite{KS} to 
show that every nontrivial surgery on $K$ is left-orderable surgery.  
Since $K$ has no cyclic period \cite[Appendix F]{Ka}, 
this example cannot be explained by the periodic construction. 
\end{enumerate}
\end{remark}

\medskip

\section{Shapes of $\mathcal{S}_{LO}(K)$ -- questions and conjectures}
\label{questions}

As we mentioned in Remark~\ref{mirror}(1), 
$0 \in \mathcal{S}_{LO}(K)$ for any knot $K$. 
If $K$ is the trivial knot then $\mathcal{S}_{LO}(K) = \{ 0 \}$, 
which has the smallest size.  
On the other hand, 
Theorems~\ref{satellite}, \ref{all rationals} and Proposition~\ref{T-pqTrs} demonstrate that there are infinitely many knots $K$ with  
$\mathcal{S}_{LO}(K) = \mathbb{Q}$, 
which has largest size. 

It seems interesting to determine the shape of $\mathcal{S}_{LO}(K)$ when it is neither $\{ 0 \}$ nor $\mathbb{Q}$. 

\begin{question}
\label{(-1,1)}
If $K$ is a nontrivial knot in $S^3$, 
then does $\mathcal{S}_{LO}(K)$ contain $(-1, 1) \cap \mathbb{Q}$? 
\end{question}

Recently Li and Roberts \cite[Corollary 1.2]{LR} prove that for any hyperbolic knot $K$, 
there exists a constant $N_K$ such that  
$\{ \frac{1}{n} \ |\ |n| > N_K \} \subset \mathcal{S}_{LO}(K)$. 

More strongly,  
we would like to ask: 

\begin{question}
\label{infty}
If $K$ is a nontrivial knot in $S^3$, 
then does $\mathcal{S}_{LO}(K)$ contain $(-\infty, 1) \cap \mathbb{Q}$ or $(-1, \infty) \cap \mathbb{Q}$? 
\end{question}

For the simplest nontrivial knot $T_{3, 2}$ (resp. $T_{-3, 2}$), 
we have $\mathcal{S}_{LO}(T_{3, 2})= (-\infty, 1) \cap \mathbb{Q}$ 
(resp. $\mathcal{S}_{LO}(T_{-3, 2})= (-1, \infty) \cap \mathbb{Q}$); 
see Remark~\ref{mirror}(2) and Example~\ref{torus knots}. 

\begin{question}
\label{trefoil}
If $\mathcal{S}_{LO}(K)= (-\infty, 1) \cap \mathbb{Q}$ 
or $\mathcal{S}_{LO}(K)= (-1, \infty) \cap \mathbb{Q}$, 
then is $K$ a trefoil knot $T_{3, 2}$ or $T_{-3, 2}$, respectively? 
\end{question}

\begin{question}
\label{maxmini}
Let $K$ be a nontrivial knot in $S^3$. 
Then does $\mathcal{S}_{LO}(K)$ have a maximum or minimum?
\end{question}

Conjecture~\ref{LspaceConjecture2} says that 
$\mathcal{S}_{L}(K)$ and $\mathcal{S}_{LO}(K)$ are complementary to each other in $\mathbb{Q}$ 
if $K$ is not a cable of a nontrivial knot. 
So let us look at the shape of $\mathcal{S}_{L}(K)$, 
which is described by Proposition 9.6 in \cite{OS4} (\cite[Lemma 2.13]{Hedden}). 

\begin{theorem}[\cite{OS4, Hedden}]
\label{L-surgery_structure}
Suppose that $K$ is a nontrivial knot and 
$\mathcal{S}_L(K) \ne \emptyset$. 
Then $\mathcal{S}_L(K) = [2g(K)-1, \infty) \cap \mathbb{Q}$ or 
$\mathcal{S}_L(K) = (-\infty, -2g(K)+1] \cap \mathbb{Q}$. 
\end{theorem}

\medskip

Theorem~\ref{L-surgery_structure} makes us expect the following explicit form of $\mathcal{S}_{LO}(K)$. 

\begin{conjecture}
\label{shape of S_LO}
Let $K$ be a nontrivial knot in $S^3$ which is not a cable of a nontrivial knot. 
Then $\mathcal{S}_{LO}(K)$ coincides with one of $\mathbb{Q}$, $(-\infty,  2g(K)-1) \cap \mathbb{Q}$ or  
$(-2g(K)+1, \infty) \cap \mathbb{Q}$. 
\end{conjecture}

Finally we give a comment on Question~\ref{trefoil} in case of $\mathcal{S}_{LO}(K)= (-\infty, 1) \cap \mathbb{Q}$; 
the other case follows by taking the mirror image. 
By the assumption $1 \not\in \mathcal{S}_{LO}(K)$. 
If Conjecture~\ref{LspaceConjecture} is true, 
then $1 \in \mathcal{S}_L(K)$ or $K(1)$ is reducible. 
The latter possibility is eliminated by \cite[Corollary 3.1] {GL},  
and hence $K(1)$ is an $L$--space. 
Then Proposition~\ref{HFPoincare} \cite[Proposition~6]{HW} below shows that $K$ is a trefoil knot $T_{3, 2}$.

\begin{proposition}[\cite{HW}]
\label{HFPoincare}
Suppose $K$ is a nontrivial knot and $K(\frac{1}{n})$ is an $L$--space.  
Then $n = 1$ $($resp. $-1$$)$ and $K$ is a trefoil knot $T_{3, 2}$ $($resp. $T_{-3, 2}$$)$. 
\end{proposition}

\bigskip

\textbf{Acknowledgements} --
We would like to thank Cameron Gordon, Hiroshi Matsuda, Yi Ni and Motoo Tange for private communications concerning 
$L$--spaces, 
and thank Joshua Greene and Adam Levine for suggesting the use of Alexander polynomials in the proof of Proposition~\ref{T-pqTrs}. 
We would also like to thank Adam Clay for informing us his curious example mentioned in Remark~\ref{infinitely many all rationals}(3), 
and Kazuhiro Ichihara for letting us know about the paper \cite{Jun}. 
Finally we would like to thank the referees for their careful reading and useful comments.

\end{document}